\documentclass[a4paper,twoside,11pt]{article}
\usepackage{pgf,tikz}
\usetikzlibrary{arrows}
\usepackage[english]{babel}
\usepackage[utf8]{inputenc}
\usepackage{amsmath}
\usepackage{amssymb}

\usepackage{amsfonts}
\usepackage{graphicx}
\usepackage{mathtools}
\usepackage[nottoc,numbib]{tocbibind}
\usepackage{tikz}
\usetikzlibrary{matrix,arrows,decorations.pathmorphing}
\usepackage{tikz-cd}
\usepackage{fancyhdr}
\usepackage{latexsym,amsmath,amssymb,amsthm, wrapfig, epsfig, dsfont}

\usepackage{pgfplots}
\pgfplotsset{width=10cm,compat=1.9}
\usepgfplotslibrary{external}

\newcommand{\R}{\mathbb{R}}
\newcommand{\C}{\mathbb{C}}
\newcommand{\Z}{\mathbb{Z}}
\newcommand\norm[1]{\left|#1\right|}

\theoremstyle{plain}
\newtheorem{propn}{Proposition}[section]
\newtheorem{thm}[propn]{Theorem}
\newtheorem{lemma}[propn]{Lemma}
\newtheorem{cor}[propn]{Corollary}

\theoremstyle{definition}
\newtheorem{defn}[propn]{Definition}

\theoremstyle{remark}

\newcommand{\bean}{\begin{eqnarray}}
\newcommand{\eean}{\end{eqnarray}}
\newcommand{\be}{\begin{displaymath}}
\newcommand{\ee}{\end{displaymath}}
\newcommand{\bea}{\begin{eqnarray*}}
\newcommand{\eea}{\end{eqnarray*}}

\numberwithin{equation}{section}
\title{Uniform  Convolution and Fourier Restriction estimates for complex polynomial curves in $\C^3$ }

\author{Conor Meade}

\date{\today}

\begin{document}
\maketitle
\thispagestyle{empty}
\begin{abstract}

\centering{}
We establish optimal $(p,q)$ ranges for two types of estimates associated to three dimensional complex polynomial curves. These are the estimates for the weighted restriction of the Fourier Transform to a complex polynomial curve, and the weighted Convolution Operator associated to a complex polynomial curve. Establishing these estimates comes down to establishing a lower bound for the Jacobian of a mapping associated to the complex curve in question. 
\end{abstract}

\thispagestyle{plain}
\setcounter{page}{1}

\lhead{Conor Meade\\
\texttt{}}
\rhead{Some three dimensional Complex Measure Convolution, and Fourier Restriction Estimates\\
\date{\today}}

\section{Introduction and Statement of Main Theorems}
Consider the complex polynomial curve $\Gamma(z):=\left(P_1(z),P_2(z),P_3(z)\right)$, for complex polynomials $P_i(z)$. If we let $L_\Gamma(z):=\left|\text{det}(\Gamma'(z),\Gamma''(z),\Gamma'''(z))\right|$, and $\lambda_\Gamma(z):=\left|L_\Gamma(z)\right|^\frac13$, we define the weighted Complex Convolution Operator associated to $\Gamma$ as: \begin{equation*}
\begin{aligned}
T_\Gamma f(\mathbf{z}):=\int_D f\left(\mathbf{z}-\Gamma(w)\right)\lambda_\Gamma(w)dw,
\end{aligned}
\end{equation*}
where $D$ is some disk of radius $r$. In this paper, we establish the following family of $(p,q)$ bounds for $T_\Gamma$:

\begin{thm}
Let $(p_\theta,q_\theta)=(\frac{6}{3+\theta},\frac{6}{2+\theta})$. Suppose $f\in L^{p_\theta}$, and $\Gamma:\C\rightarrow \C^3$ is a degree $N$ polynomial curve. Then there exists a constant $C=C(N,\theta)$ such that: $$||T_\Gamma f||_{L^{q_\theta}(\C^3)}\leq C||f||_{L^{p_\theta}(\C^3)}, $$
for all $\theta \in (0,1)$.

Furthermore, $T_\Gamma$ is of restricted weak-type at $(p_0,q_0)$, and at $(p_1,q_1)$, i.e. for measurable subsets of $\C^3$, $E$ and $F$, there exists $C=C(N)$ such that $$\langle T\chi_E,\chi_F\rangle\leq C\norm{E}^{\frac 1p_{i}}\norm{F}^{1-\frac{1}{q_{i}}},$$ for $i=0$ or $1$.
\end{thm}
Utilising this theorem, some easy to obtain estimates, and the Marcinkiewicz Interpolation Theorem, we can obtain a larger range of local estimates.
\newpage
\begin{cor}

With $(p_\theta,q_\theta)$ as in Theorem 1.1, let $\mathcal{R}$ be the convex hull of the points $\left(1,1\right)$, $\left(0,0\right)$, $\left(p_0^{-1},q_0^{-1}\right)$, and $\left(p_1^{-1} ,q_1^{-1}\right)$. 

If $(p^{-1},q^{-1})\in\mathcal{R}\setminus\left\{\left(p_0^{-1},q_0^{-1}\right),\left(p_1^{-1} ,q_1^{-1}\right)\right\}$, then for all complex polynomial curves, $\Gamma$, there exists a constant $C=C(\Gamma,r,p,q)$, where $r$ is the radius of the disk associated to $T_\Gamma$, such that for all functions $f\in L^p(\C^3)$: $$
||{ T_\Gamma f }||_{L^q(\C)}\leq C||{f}||_{L^p(\C)}.$$

\end{cor}

Much work has been done on estimates of this sort are known when $T_\Gamma$ is replaced with it's real analogue. In the real case, in arbitrary dimension, $d$, $\Gamma$ is a real polynomial curve, the integration is over some interval, and the weight function $\lambda_\Gamma$ is replaced with $\lambda_\Gamma(t):=L_\Gamma^\frac{2}{d(d+1)}(t)$. Furthermore, our constant $C$ must depend on the dimension we are in for these estimates to hold.

Usually, these estimates have been for a particular family of polynomial curves, or in a particular dimension. For example, in 1998, in \cite{Christ1} one of the early papers on this topic, Christ established this result for any moment curve in any dimension. That is to say he established it for $\Gamma(t)=\left(t,t^2\cdots t^n\right)$ in $\R^n$. It should be noted that for curves like this, we have $\lambda_\Gamma(t)=n!$, which simplifies the problem. This can be considered the most non-degenerate case, as all derivatives of $\Gamma$ are linearly independent. Later papers introduced the weight $\lambda_\Gamma$ in order to deal with problematic singular points of $\Gamma$, where the derivatives become linearly dependent. At such points, $\lambda_\Gamma$ vanishes, allowing us to mitigate the effect of such points during the integration. 

Further estimates were established for other families of curves, but it was in 2002, in \cite{Oberlin}, that Oberlin established the real estimate for general polynomial curves in two dimensions, and in three dimensions, established it for curves of the form $\gamma(t)=\left(t,P_1(t),P_2(t)\right)$. In \cite{DendrinosLaghiWright}, this three dimensional case was extended to general polynomial curves, and in \cite{Stovall2}, this was extended to general higher dimensions.

The first such estimates for surfaces instead of curves can be found in \cite{DruryGuo},     in which Drury and Guo proved a class of estimates of this type, but instead of convolving with a measure arising from a polynomial curve, their estimates were pertaining to the case where $\Gamma$ was instead a two dimensional surface in $\R^4$, of the form $(x,y,\phi_1(x,y),\phi_2(x,y))$, with certain conditions on $\phi_1$ and $\phi_2$. 

In \cite{ChungHam}, Chung and Ham established similar results, considering surfaces that could be represented as complex curves in $\C^d$. They established these bounds for curves of the form $\left(z,z^2\cdots z^N\right)$ in $\C^d$, and $\left(z,z^2,\phi(z)\right)$ in $\C^3$, for analytic $\phi$. They do this by using a more powerful analogue to our Lemma 3.1, established in the complex Fourier Restriction paper \cite{BakHam}, involving a lower bound in terms of the arithmetic means of $\{L_\Gamma(z_i)\}$, as opposed to the geometric mean. This method suffices for the complex curves in question in \cite{ChungHam}, but in general does not hold. 
In this paper, we do not utilise a lower bound in terms of the arithmetic mean, but instead adapt the geometric mean bound deriving from \cite{DendrinosWright} into the complex case, and achieve the overall estimate using that.

Note that while the two-dimensional estimate has not been established for general complex polynomial curves, these estimates can be derived from the procedure described in our proof of Theorem 1.1, with simple modifications. The explicit description of this is omitted to aid the presentation of the three dimensional case.

The second types of estimates we wish to establish are pertaining to the restriction of the Fourier transform to $\Gamma$, and are formulated as follows:
\begin{thm}

Let $\mathcal{L}$ be the line segment joining the points $(1,0)$ and $(\frac67,\frac 67)$, excluding the endpoint  $(\frac67,\frac 67)$. If $(p^{-1},q^{-1})\in\mathcal{L}$, then for all complex polynomial curves, $\Gamma$, of degreee $N$, there exists a constant $C=C(p,N)$, such that for all functions $f\in L^q$:
$$||\hat{f}(\Gamma(z))||_{L^q(\lambda_\Gamma(z) dz)}\leq C||f||_{L^p(dz)}, $$ 
where $\hat{f}$ denotes the Fourier Transform of $f$.
\end{thm}

Here, we should note that $\hat{f}(\Gamma(z))$ is understood as $\hat{f}(\Gamma(x,y))$, where $z=x+iy$, and we regard $\Gamma$ as a real surface, in order to ensure a sensible notion of convergence for it's Fourier transform. This detail will be omitted, unless in the case of potential confusion.

Note that in \cite{BakHam}, Bak and Ham established the optimality of the weight function $\lambda_\Gamma$ in the complex case, in the sense that it is the largest possible weight such that we can obtain this range of $(p,q)$ values. Furthermore, in \cite{WrightUnpublished}, we see this this is the largest possible $(p,q)$ range for estimates of this type. 

With regard to the history of this problem, initial investigations pertained to the real case, with $\gamma$ being a real curve instead of a surface . In \cite{Sjoln}, Sj\"oln demonstrated that the difficulty in acquiring such bounds was in large part due to potential oscillations of the quantity $L_\Gamma$. Considering the curve $\Gamma(t)=(t,e^{-\frac 1t}\sin{\frac 1t})$ for $t$ close to $0$ demonstrates this issue. In this paper, Sj\"oln avoided issues like this by restricting $\Gamma$ to convex curves, ensuring the single signedness of $L_\Gamma$. 

There are several results that could be cited here, such as \cite{BakOberlinSeeger}, which establishes the estimate for $\gamma(t)$ being an abitrary monomial curve in $\R^d$, that is to say each of the components of $\Gamma$ is some monomial. Furthermore, in the same paper, these estimates are established for $\Gamma(t)$, being of what is called ``simple type'', which is to $\Gamma(t)=\left(t,t^2\cdots P(t)\right)$ for polynomial $P$.

Usually, attempts to obtain these restriction estimates aim for one of two different ranges, which each depend on the dimenson. There is the full range, in which we have $(p,q)$ satisfying $p\in [1,\frac{d^2+d+2}{d^2+d})$, and $p'=\frac{d^2+d+2}{2}q$, and the shorter range where we still have $p'=\frac{d^2+d+2}{2}q$, but we only have $p\in [1,\frac{d^2+d+2}{d^2+2d-2})$. Note that here, $p'$ denotes the Hölders conjugate of $p$. In \cite{DendrinosWright}, the short range was established for arbitrary real polynomial curves, in any dimension, being done by decomposing $\R$ into intervals on which the quantity $L_\gamma$ behaved as a monomial on any interval arising from this decomposition, and this method was utilised and improved in \cite{Stovall} to cover the full range, giving a complete answer to this type of estimate in the real case, in the context of boundedness between $L^p$ spaces.

The complex case also has some results established. In \cite{Oberlin2} Oberlin established these bounds over the full range for surfaces of the form $\left(z,\phi(z)\right)$ with some non-degeneracy conditions on $\phi$. 

In \cite{BakHam}, this estimate was established on the full range for complex curves of the form $\Gamma(z)=\left(z,z^2\cdots z^N\right)$ in $\C^d$, for $d\geq 3$, and for $\Gamma(z)=\left(z,z^2,\phi(z)\right)$ in $\C^3$, for arbitrary analytic $\phi$.

Before moving onto the process of proving Theorem 1.1 and 1.3, we first introduce some commonly accepted notation for these kinds of problems. We say that, for positive quantities $A$ and $B$, we have $A\lesssim B$ if $A\leq C B$, for some constant $C$. Here, $C$, the implicit constant, can depend on whatever parameters are appropriate for the problem in question.

As Theorem 1.1 follows from a single restricted weak-type estimate at $(p_0,q_0)$, when proving this theorem, we only want our implicit $C$ to be dependant on $N$, the degree of $\Gamma$. In Theorem 1.3, we prove the estimate via many iterative estimates, so we would expect $C$ to depend on $p$ in this context. We also say $A\sim B$ if $A\lesssim B\lesssim A$. 

\section{Decomposition of $\C$}
To achieve the desired results, we first partition $\C$ into various convex sets. This is done by utilising two decompositions as described in \cite{DendrinosWright}. In \cite{DendrinosWright}, these decompositions are described in relation to $\R$, but can be extended to $\C$. The first of these decompositions will be referred to as D1.

D1: Given a polynomial $Q(z)$, we can decompose $\C$ into a bounded number of convex sets, $B_i$. On each of these $B_i$, we have $\norm{Q(z)}\sim c(Q)\norm{z-b_i}^{k_i}$, for $b_i$ some root of $Q(z)$. We will frequently omit the subscript of $k$ and $b$ for the sake of convenience.

To see how this decomposition functions, we first write 
${\displaystyle Q(z)=A\prod_{j=1}^{d'}(z-\eta_j)^{\alpha_j}}$, with $d'$ being the amount of distinct roots of $Q$. Here, these roots $\eta_j$ are written such that $\norm{\eta_1}\leq\norm{\eta_2}\leq\cdots\leq\norm{\eta_{d'}}$. We now decompose $\C$ into a union of $d'$ convex sets, given by the Voronoi diagram associated to $\{\eta_1,\eta_2,\cdots,\eta_{d'}\}$. On these sets, we will have our complex centers $b$ be the root associated to this set. Formally these sets are denoted $S(\eta_j)=\{z\in\C| \norm{z-\eta_j}\leq\norm{z-\eta_i}\text{ for all }i\neq j \}$. 

Each of these sets are further divided into finitely many sectors, centered at $b$. For $n\geq 1$ denote $\Delta^n_\varepsilon(\eta_j)=\{z-\eta_j=re^{i\theta}\in S(\eta_j)|(n-1)\varepsilon\leq\theta\leq n\varepsilon\}$. Here, $\varepsilon$ is chosen such that $\frac{2\pi}{\varepsilon}$ is an integer. Specifically we need to take $\varepsilon\leq\frac{\pi}{8}$, which is necessary to decompose the intersection of some annuli and these sectors into convex sets. While this is sufficient for $D1$ to function as described, we will see that we also need to choose $\varepsilon$ to be dependent on $d$, the degree of our polynomial curve. Thus the choice of $\varepsilon$ is currently not specified. 

Fixing one such $\Delta^n_{\varepsilon}(\eta_j)=\Delta(\eta_j)$, we relabel $\eta_j=\hat{\eta}_1$, and the remaining roots as $\hat{\eta}_k$ so that we have $\norm{\hat{\eta}_1-\hat{\eta}_2}\leq\norm{\hat{\eta}_1-\hat{\eta}_3}\leq\cdots\leq \norm{\hat{\eta}_1-\hat{\eta}_{d'}}$. That is to say we order them according to their distance from $\eta_j$ We let, for $i\geq 2$, $T^n_{i,j}=\{z\in\Delta(\eta_j)|\norm{z-\hat{\eta}_1}\leq\frac{1}{2}\norm{\hat{\eta}_1-\hat{\eta}_i}\}$. 
Labelling $T_{1,j}=\varnothing$, and $T_{d'+1,j}=\Delta(\eta_j)$, we see that $T_{1.j}\subseteq T_{2.j}\subseteq\cdots \subseteq T_{d'.j}\subseteq T_{d'+1.j}=\Delta(\eta_j)$. Also note that if $z\in T^n_{i,j}$, then $\norm{\hat{\eta}_1-\hat{\eta}_{k}}\sim\norm{z-\hat{\eta}_{k}}$ for all $k\geq i$, and if $z\notin T^n_{i,j}$, then $\norm{z-\hat{\eta}_1}\sim\norm{z-\hat{\eta}_k} $ for all $k\leq i$.

Letting, for $1\leq i\leq d'$, $I^n_{i,j}=T^n_{i+1,j}\backslash T^n_{i,j}$, we can see that $\Delta^n_\varepsilon(\eta_j)=\bigcup_{k=1}^{d'}I^n_{k,j}$, and on any fixed $I_{k,j}$, we have:
$$
\norm{Q(z)}=\norm{A}\prod_{l=1}^{d'}\norm{z-\eta_l}^{\alpha_l}\sim \norm{z-\hat{\eta}_1}^{\hat{\alpha}_1+\hat{\alpha}_2+\cdots +\hat{\alpha}_k}\norm{A}\prod_{l=k+1}^{d'}\norm{\hat{\eta}_1-\hat{\eta}_l}^{\hat{\alpha_l}},
$$

for $\hat{\alpha}_j$ being the multiplicity of $\hat{\eta}_j$.

So indeed, writing 
${\displaystyle\C=\bigcup_{l=1}^{\frac{2\pi}{\varepsilon}d'^2} I_l=\bigcup_{j=1}^{d'}\bigcup_{n=1}^{\frac{2\pi}{\varepsilon}} \bigcup_{i=1}^{d'} I^n_{i,j}}$, we have on $B_l=I^n_{i,j}$ that:
\begin{equation}
\begin{aligned}
\norm{Q(z)}\sim c(Q)\norm{z-b_l}^{k_l},
\end{aligned}
\end{equation}
where $b_l=b_{i,j,n}=\eta_j$, and $k_l=k_{i,j,n}=\hat{\alpha}_1+\hat{\alpha}_2\cdots +\hat{\alpha}_i$.

Note that these sets are intersections of some convex sets, and annuli centered at $b_l$, and therefore are not necessarily convex. This can be remedied by ``convexifying'' our annuli. Note that we have $\varepsilon<\frac{\pi}{4}$, and we can scale down the inner radius and scaling up the outer radius of each annuli. As long as these scalings are given by a constant multiplicative factor, $B(\varepsilon)$, the new sets retain the relevant properties of our original sets, but are no longer disjoint from their 'neighbouring' annuli. This means that the estimate (2.1) for each of the original annuli will hold on the intersection of the two scaled annuli.

Consider the bisector of $\Delta^n_\varepsilon(\eta_j)$, and its intersections with the innermost boundary of the intersection of these annuli. Note that this boundary arises from the annulus that comes \emph{later} in our sequence of annuli, because of our thickening. By choosing $B$ suitably, we can control the thickness of the two annulis' intersection, and ensure that the tangent to the innermost boundary divides the intersection into two sets,  without intersecting the outer boundary. This then allows us to replace the our annuli with convex sets by replacing the curved boundaries of these sets with these tangent lines.

For the next decomposition, referred to as D2, we require a polynomial $Q(z)$, and a centre, which is some complex number $b$.

D2: Given a polynomial $Q(z)$, and complex number $b$, and some convex set $J$, then $J$ can be decomposed into disjoint convex sets. These convex sets are of two types, either dyadic or gap. On Dyadic sets, $\norm{z-b}\sim c_1(Q)$, and on Gap sets, $Q(z)\sim c_2(Q)\norm{z-b}^k$, for some positive integer $k$. Also, for $Q(z+b)=\Sigma c_iz^i$, if $c_j=0$, then there are no gap intervals such that $\norm{Q(z)}\sim c_2(Q) \norm{z-b}^j$. 

Note that, if we label the collection of Dyadic sets as $\{D_j\}$, and the gaps as $\{G_j\}$we have:\begin{itemize}
    \item $D_j=\{z\in\C|A^{-1}\norm{z_j}<\norm{z}<A\norm{z_j}\}$,
    \item $G_j=\{z\in\C|A\norm{z_j}<\norm{z}<A^{-1}\norm{z_{j+1}}\}$,
 \end{itemize}
 for $A\in\R$, and $z_j$ being the roots of $Q(z)$, ordered by magnitude.

D2 originates from \cite{CarberryRicciWright}, and describes a decomposition of $\R$. However, if we  restrict ourselves to sectors like $\Delta^n_\varepsilon(\eta_j)$, where $\eta_j$ are the roots of $Q$, and Voronoi sets $S(\eta_j)$, then the method can be extended to the complex case. Note again however that this will give us disjoint annuli in the complex case, so we repeat the thickening of these annuli, and converting them into corresponding convex sets as in D1. 

For now, we will assume that $C(Q)=c_1(Q)=c_2(Q)=1$, with the justification for this assumption being apparent whe we introduce the main inequality we seek to establish.

We also make use of the integral representation of the Jacobian in \cite{DendrinosWright}. This states that, in 3-dimensions, for $L_3(z)=L_\Gamma(z)$, $L_2(x)=\text{det}\left(\begin{array}{cc}
    P_1'(z) & P_1''(z) \\
     P_2'(z)& P_2''(z)
\end{array}\right)$, and $L_1(z)=P_1'(z)$, then if we consider the mapping $\Phi_\Gamma(z_1,z_2,z_3)=\Gamma(z_1)+\Gamma(z_2)+\Gamma(z_3)$, we have:

\begin{equation}
\begin{aligned}
J_{\Phi_\Gamma}(z_1,z_2,z_3)=\prod_{s=1}^3L_1(z_s)\int_{z_1}^{z_2}\int_{z_2}^{z_3}\prod_{s=1}^2\frac{L_2(w_s)}{L_1(w_s)^2}\int^{w_2}_{w_1}\frac{L_1(y)L_3(y)}{L_2(y)^2}dydw_2dw_1,
\end{aligned}
\end{equation}

Where $J_{\Phi_\Gamma}$ is the real jacobian of the mapping $\Phi_\Gamma$, and the complex bounds of integration, $u,v$ in each of the above integrals are to be interpreted as the line integral between those two complex numbers along the curve given by $\Gamma(t)=ut+(v-u)t$ for $t\in[0,1]$.

From here, we use D1 and D2 to decompose $\C$ into convex sets so that $L_1$, $L_2$, and $L_3$ can be treated as either monomials or constants. This is done as follows.

First, we use D1 with respect to $L_3$, so that we get a finite amount of convex sets such that $\norm{L_3} \sim \norm{z-b}^k$ on each of these sets, with $k$ and $b$ depending on the sets. Next, we divide this family of sets into two types of sets, referred to as $T_0$ or $T_1$ sets. To do this, we first apply D2 to each set already obtained with respect to our center from the first decomposition, and the polynomial $L_1$.

$T_0$ sets are those sets are gap annuli, on which $\norm{L_1(z)}\sim\norm{z-b}^{k_0}$, and $\norm{L_3(z)}\sim\norm{z-b}^k$.

$T_1$ sets are dyadic annuli those on which $\norm{z-b}\sim 1$, so that $\norm{L_3(z)}\sim 1$. On these sets we apply D1 with respect to $L_1$ to decompose these sets into a family of convex sets on which $\norm{L_1(z)}\sim\norm{z-b'}^{k_1}$.

Finally, we apply D2 again to the $T_0$ and $T_1$ sets, but this time, with respect to $L_2$, and the center $b$, and $b'$ respectively. This divides the $T_0$ sets into two new families of sets, referred to as $T_{00}$ sets and $T_{01}$ sets.

$T_{00}$ sets are those on which $\norm{L_2(z)}\sim\norm{z-b}^{k_{00}}$, while we still have $\norm{L_1(z)}\sim\norm{z-b}^{k_0}$, and $\norm{L_3(z)}\sim\norm{z-b}^k$

$T_{01}$ sets are those on which $\norm{z-b}\sim 1$, so that $\norm{L_1(z)}\sim \norm{L_3(z)}\sim 1$. Like previously, we apply D1 with respect to $L_2$ to decompose these sets into a family of convex sets on which $\norm{L_2(z)}\sim\norm{z-b''}^{k_{01}}$.

We decompose the $T_1$ sets precisely the same way to get families of sets named $T_{10}$ and $T_{11}$.

On $T_{10}$ sets, we have $\norm{L_3}\sim 1$, $\norm{L_1}\sim\norm{z-b'}^{k_1}$
and $\norm{L_2(z)}\sim\norm{z-b'}^{k_{10}}$. 

On $T_{11}$ sets, we have $\norm{L_1(z)}\sim \norm{L_3(z)} \sim 1$, and $\norm{L_2(z)}\sim\norm{z-b''}^{k_{11}}$.

Note that the various $b$ and $k$ values depend on the specific set within these families that they belong to, as opposed to every set of type $T_{01}$ having the same $k_{01}$ value for example.

Furthermore, we shall write $b'$ and $b''$ as just $b$ for simplicity for the rest of this discussion. Note to make use of these estimates, we need the following lemma.

\begin{lemma}
For all $z_1,z_2,z_3$ in one of the sets in our decomposition, we have that $$|J_{\Phi_\Gamma}(z_1,z_2,z_3)|\sim\prod_{s=1}^3|L_1(z_s)|\left|\int_{z_1}^{z_2}\norm{\int_{z_2}^{z_3}\prod_{s=1}^2\frac{|L_2(w_s)|}{|L_1(w_s)|^2}\left|\int^{w_2}_{w_1}\frac{|L_1(y)||L_3(y)|}{|L_2(y)|^2}dy\right|dw_2}dw_1\right|.$$
\end{lemma}

\emph{Proof:} To see this, we introduce the notion of a sector contained function. 

\begin{defn}
An  $\varepsilon$-sector contained function of $B$ is a function, $f$, such that for all $z\in B$, we can  write $f(z)=r(z)e^{i\theta(z)}$, with $\theta(z)\in[\theta_0,\theta_0+\varepsilon]$. That is to say, the image of $B$ under $f$ is contained within a sector of aperture $\varepsilon$.
\end{defn}

We will omit the reference to the set $B$ when we refer to function that satisfy this definition, understanding that a $\varepsilon$-sector contained function refers to a $\varepsilon$-sector contained function of $D$, where $D$ is any set in our decomposition.

Note that the reason we introduce the notion of these functions is because for $f$ a $\varepsilon$-sector  contained function, we have $\norm{\int_u^v f(z)dz} \sim\norm{\int_u^v \norm{f(z)}dz}$ for $u,v$ in any set of our decomposition, for $\varepsilon=\varepsilon(d,N)<\frac \pi 2$.

To see this, note that we immediately have that $\norm{\int_u^v f(z)dz}\leq\norm{\int_u^v \norm{f(z)}dz}$, as $$\norm{\int_u^v f(z)dz}=\norm{(v-u){\int_0^1 f(z(t))dt}}\leq\norm{(v-u){\int_0^1\norm{ f(z(t))}dt}}=\norm{\int_u^v \norm{f(z)}dz}.$$

To establish the second inequality, notice that by factoring out some complex number of unit length, we can assume the sector $f(z)$ is contained in a sector which has one ray in the direction of the positive real axis, and the second ray being being given by is the line containing the complex numbers with argument $\varepsilon$.

We therefore have:
\begin{equation*} 
\begin{split}
\norm{\int_u^v f(z)dz} & = \norm{(v-u)\int_0^1 \text{Re}(f(z(t)))+i\text{Im}(f(z(t)))dt}\\
 & \geq\norm{(v-u)\int_0^1 \text{Re}(f(z(t)))dt}\\
 &\sim\norm{\int_u^v\norm{f(z)}dz}.
\end{split}
\end{equation*}
Therefore, to establish our lemma, we would like to show that each integrand is an $\varepsilon$-sector contained function. However, because the argument of $w_2-w_1$ cannot be controlled, this is not quite true, so we instead show that the integrand is the product of some $\varepsilon$ contained function, and $(w_2-w_1)$, which we will see is sufficient to establish our lemma.

In \cite{ChungHam}, it was shown in step 2 of lemma 4.2 that for a polynomial \newline$P(z)=\prod_{j=1}^d (z-\eta_j)$, if we are in $\Delta(\eta_j)_\varepsilon\bigcap S(\eta_j)$, as in the terminology of D1, then we have, by factoring $P(z)=g_j(z)(-1)^{d-j}z^j\prod_{l=j+1}^d\eta_l$, that on the $j$-th gap or dyadic annulus, we can decompose these sets into a bounded number of sets so that $g_j(z)$ is contained in a sector of aperture less that $\varepsilon$. As $(-1)^{d-j}z^j\prod_{l=j+1}^d\eta_l$ will be contained in a sector of aperture $j\varepsilon$, we have that $P(z)$ will be in a sector of aperture at most $(d+1)\varepsilon$ for $z$ in any of our sets. Notice that both decomposition techniques $D1$ and $D2$ involve restricting our sets to  $\Delta(\eta_j)_\varepsilon\bigcap S(\eta_j)$, and on each of our sets we have applied $D1$ or $D2$ to each of $L_1,L_2$ and $L_3$. Therefore, we have no concern about utilising this result for $L_1,L_2$ and $L_3$ simultaneously.

So therefore, for any set in our decomposition, we can decompose it further to ensure each of our $L_i$ are $(d_i+1)\varepsilon$-sector contained functions, for $d_i$ the degree of $L_i$. Therefore, we have that the rational function $\frac{L_1L_3}{L_2^2}$ is a $4(d+1)\varepsilon$-sector contained function, and $\frac{L_2}{L_1^2}$ is a $3(d+1)\varepsilon$-sector contained function where $d$ is the maximum of the $d_i$ values.

Lastly, we wish to observe that $\int^{w_2}_{w_1}f(y)dy$ is the product of a $\varepsilon$-sector contained function and $(w_2-w_1)$ if $f$ is a $\varepsilon$-sector contained function, and $w_1,w_2$ are contained within one of our decomposed sets. This follows from the convexity, and closedness of sectors.

\begingroup\fontsize{9.5}{12}
Letting $\int^{w_2}_{w_1}f(y)dy=(w_2-w_1)\text{lim}_{n\rightarrow\infty}\frac 1n\Sigma_{j=1}^nf(y(\frac jn))$, for $y(t)=w_1+(w_2-w_1)t$. \newline Now, as $f(y)$ is in a sector for all $y$ on this line segment, as the sector is convex, each $\frac 1n\Sigma_{j=1}^nf(y(\frac jn))$ is in the sector, so therefore the limit of this sequence is in the sector. 

Therefore, we have that $\int^{w_2}_{w_1}f(y)dy=(w_2-w_1)F(z)$ , where $F(z)=\int_0^1(f(y(t))dt$ is a $\varepsilon$-sector contained function. Now note that we can write:
\begin{equation*}
\begin{split}
J_{\Phi_\Gamma}(z_1,z_2,z_3)&=\prod_{s=1}^3L_1(z_s)\int_{z_1}^{z_2}\int_{z_2}^{z_3}\prod_{s=1}^2\frac{L_2(w_s)}{L_1(w_s)^2}\int^{w_2}_{w_1}\frac{L_1(y)L_3(y)}{L_2(y)^2}dydw_2dw_1\\
 & =\left(\prod_{s=1}^3L_1(z_s)\right)(z_2-z_1)(z_3-z_2) I\\
 \text{where:}\\
 I&:=\int_{0}^{1}\int_{0}^{1}\prod_{s=1}^2\frac{L_2(w_s)}{L_1(w_s)^2}(w_2-w_1))\int^{1}_{0}\frac{L_1(y)L_3(y)}{L_2(y)^2}dt_3dt_2dt_1.
\end{split}
\end{equation*}
Here, note that $w_1=w_1(t_1)$, $w_2=w_2(t_2)$, and $y=y(t_1,t_2,t_3)$.

We now wish to bound the size of the remaining integral $I$. Let $\alpha(t_1,t_2)$ and $\beta(t_1,t_2)$ be real valued functions such that $w_2-w_1=\alpha+i\beta$. Furthermore, as these polynomials, after our decomposition, are section-contained functions, and the integral along the line from $w_1$ to $w_2$ of a sector contained function is a sector contained function multiplied by $(w_2-w_1)$. Therefore, the above integral can be written as: 

\begin{equation*}I=\int_{0}^{1}\int_{0}^{1}(\alpha+i\beta)g(t_1,t_2)dt_2dt_1,
\end{equation*}
where $g$ is a 
$(7(d+1))\varepsilon$-sector contained function. That is to say $g(t_1,t_2)=\xi(t_1,t_2)+i\eta(t_1,t_2)$, for real functions $\xi$ and $\eta$ with $|g(t_1,t_2)|\sim\xi(t_1,t_2)$. Note that we make the following assumptions about the points $z_1,z_2,z_3$. These assumptions can be made without loss of generality via renaming these points for the first two, and factoring out a unit complex number for the last. First, we assume that in the triangle $z_1,z_2,z_3$, the angle at $z_2$, which we label $\theta$, is the largest in the triangle. Next, we assume $|z_3-z_2|\geq|z_2-z_1|$. Lastly we assume $z_2-z_1$ is a positive real number. \endgroup

Note that this immediately gives that $\beta$ is single signed, and $\theta\in[\dfrac{\pi}{3},\pi]$. We compute our estimate by splitting into cases when $\theta$ is acute, and when it is obtuse.

Case i): If $\theta\in[ \dfrac{\pi}{3},\dfrac{\pi}{2}]$:

Note in this section, we will assume $\beta\geq0$. The case where $\beta<0$ is dealt with by taking $|\beta|$ instead of $\beta$ when applying the triangle inequality in the following argument:
\begin{equation*}
\fontsize{10pt}{12pt}
\begin{split}
|I|&\geq\norm{\text{Im}\left(\int_{0}^{1}\int_{0}^{1}(\alpha+i\beta)g(t_1,t_2)dt_2dt_1\right)}\\
& \geq\norm{\iint_{\beta>\frac{|\alpha|}{2}}(\beta\xi+\eta\alpha)dt_2dt_1}\\
 &\geq\norm{\iint_{\beta>\frac{|\alpha|}{2}}(\beta\xi-|\eta||\alpha|)dt_2dt_1}\\
 &\sim\iint_{\beta>\frac{|\alpha|}{2}}\beta\xi dt_2dt_1.
\end{split}
\end{equation*}
Now, on the set we restrict our integral to, $\beta\sim |w_2-w_1|$, and we have 
\begin{equation*}
\begin{split}
\xi\sim |g(t_1,t_2)|&=\norm{\prod_{s=1}^2\frac{L_2(w_s(t_s))}{L_1(w_s(t_s))^2}\int^{1}_{0}\frac{L_1(y(t_3))L_3(y(t_3))}{L_2(y(t_3))^2}dt_3}\\
 & =\prod_{s=1}^2\frac{|L_2(w_s(t_s))|}{|L_1(w_s(t_s))^2|}\norm{\int^{1}_{0}\frac{L_1(y(t_3))L_3(y(t_3))}{L_2(y(t_3))^2}dt_3}.\\
\end{split}
\end{equation*}

Now, as the integrand here is a sector contained function, we can move the absolute value signs inside the integral, yielding
\begin{equation}
\begin{split}
|g(t_1,t_2)|&\sim\prod_{s=1}^2\frac{|L_2(w_s(t_s))|}{|L_1(w_s(t_s))^2|}\left|\int^{1}_{0}\frac{|L_1(y(t_3))||L_3(y(t_3))|}{|L_2(y(t_3))^2|}dt_3\right|\\
 & \sim\prod_{s=1}^2|w_s(t_s)|^{\sigma_2}\left|\int^{1}_{0}|y(t_3)|^{\sigma_3}dt_3\right|.
\end{split}
\end{equation}

Subbing this into $(2.4)$ gives $$|I|\gtrsim \left|\int_0^1\int_{t_2(t_1)}^1|w_2-w_1|\prod_{s=1}^2|w_s(t_s)|^{\sigma_2}\left|\int_0^1|y(t_3)|^{\sigma_3}dt_3\right|dt_2dt_1\right|,$$

where $t_2(t_1)$ is the smallest value of $t_2$ for which the complex number $w_2-w_1$ has imaginary part greater than half the absolute value of it's real part. That is to say, $t_2(t_1)$ is the smallest value such that $\beta>\dfrac{|\alpha|}{2}$. Note that by our assumptions on $\theta$ and the relative lengths of our line segments, we have that $t_2(t_1)<\frac12$ for all $t_1$.

From here we will proceed to estimate the $t_1$ integrand, by considering 
\begin{equation*}
\begin{split}
G(t_1)&=\int_{t_2(t_1)}^1|w_2-w_1|\prod_{s=1}^2|w_s(t_s)|^{\sigma_2}\left|\int_0^1|y(t_3)|^{\sigma_3}dt_3\right|dt_2\\
 & =\int_{2t_2(t_1)}^1|w_2-w_1|\prod_{s=1}^2|w_s(t_s)|^{\sigma_2}\left|\int_0^1|y(t_3)|^{\sigma_3}dt_3\right|dt_2\\
 &\qquad+\int_{t_2(t_1)}^{2t_2(t_1)}|w_2-w_1|\prod_{s=1}^2|w_s(t_s)|^{\sigma_2}\left|\int_0^1|y(t_3)|^{\sigma_3}dt_3\right|dt_2\\
 &:=G_1+G_2.
\end{split}
\end{equation*}

Now note that $L=z_2-w_1$ is a positive real number, and for all the $t_2$ in the $G_2$ integral, we have $L\leq|w_2-w_1|\leq2L$, and also that this relation holds for $t_2<t_2(t_1)$. Furthermore, we also note that the triangle with vertices $z_1,z_2,w_2(2t_2(0))$ is contained in the circle centered at $z_2$, with radius $\varepsilon|z_2|$. This can be seen by considering the case with $z_1,z_2$ are positive real numbers,  $z_3$ is on the other ray of our sector, $\theta=\frac\pi 2$, and $|z_2-z_1|=|z_3-z_2|$, as this arrangement maximises the distance from $z_2$ and $w_2(2t_2(0))$. Note for all points $z$ in this circle, we have $|z|\sim|z_2|$. We therefore have:
\begin{equation*}
\begin{split}
G_2&\sim\int_{t_2(t_1)}^{2t_2(t_1)}|w_2-w_1||z_2|^{2\sigma_2+\sigma_3}\left|\int_0^1dt_3\right|dt_2\\
 & \sim\int^{t_2(t_1)}_0|w_2-w_1||z_2|^{2\sigma_2+\sigma_3}\left|\int_0^1dt_3\right|dt_2\\
 &\sim\int^{t_2(t_1)}_0|w_2-w_1|\prod_{s=1}^2|w_s(t_s)|^{\sigma_2}\left|\int_0^1|y(t_3)|^{\sigma_3}dt_3\right|dt_2\\
 &:=G_3.
\end{split}
\end{equation*}

Thus, we can write:
\begin{equation*}
\begin{split}
G(t_1)&\sim G_1+G_2+G_3\\
 & =\int_0^1|w_2-w_1|\prod_{s=1}^2|w_s(t_s)|^{\sigma_2}\left|\int_0^1|y(t_3)|^{\sigma_3}dt_3\right|dt_2.
\end{split}
\end{equation*}

Integrating both sides with respect to $t_1$ as $t_1$ goes from $0$ to $1$ yields $$|I|\gtrsim\norm{\int_0^1\int_0^1|w_2(t_2)-w_1(t_1)|\prod_{s=1}^2|w_s(t_s)|^{\sigma_2}\left|\int_0^1|y(t_3)|^{\sigma_3}dt_3\right|dt_2dt_1},$$

which is our desired bound.

Case ii) $\theta\in[\frac{\pi}{2},\pi]$.

This case is simpler, as now we always have $\alpha>0$. Using the same notation as earlier, we have 
\begin{equation} \label{eq1}
\begin{split}
|I|&=\norm{\int_0^1\int_0^1(\alpha+i\beta)(\xi+i\eta)dt_2dt_1}\\
 & =\norm{\int_0^1\int_0^1\alpha\xi-\beta\eta+i(\beta\xi+\alpha\eta)dt_2dt_1}.
\end{split}
\end{equation}

Note  that we also have, as $\beta$ is single signed, and $\alpha$ is positive, $$
\int_0^1\int_0^1\alpha\xi dt_2dt_1+
\left|\int_0^1\int_0^1\beta\xi dt_2dt_1\right|=\left|\int_0^1\int_0^1(\alpha+|\beta|)\xi dt_2dt_1\right|,$$ $$\sim\left |\int_0^1\int_0^1|w_2(t_2)-w_1(t_1)||g(t_1,t_2)|\xi dt_2dt_1\right|:=G.$$

If $\left|\int_0^1\int_0^1\beta\xi dt_2dt_1\right|\geq\int_0^1\int_0^1\alpha\xi dt_2dt_1$, then we have: $$|\text{Im}(I)|=\left|\int_0^1\int_0^1\beta\xi+\alpha\eta dt_2dt_1\right|\geq\left|\int_0^1\int_0^1\beta\xi dt_2dt_1\right|-\varepsilon\int_0^1\int_0^1\alpha\xi dt_2dt_1$$
$$\sim \left|\int_0^1\int_0^1\beta\xi dt_2dt_1\right|\gtrsim G .$$

While if $\int_0^1\int_0^1\alpha\xi dt_2dt_1\geq\left|\int_0^1\int_0^1\beta\xi dt_2dt_1\right|$, then by the same line of reasoning we have $$|\text{Re}(I)|=\left|\int_0^1\int_0^1\alpha\xi-\beta\eta dt_2dt_1\right|\gtrsim G.$$

So indeed this gives us in either case that
$$|I|\gtrsim|\int_0^1\int_0^1|w_2-w_1|\prod_{s=1}^2|w_s(t_s)|^{\sigma_2}\left|\int_0^1|y(t_3)|^{\sigma_3}dt_3\right|dt_2dt_1|.$$
as again can move the absolute value around the inner integral into it, as the integrand the a sector contained function.

This completes the lemma, as  in all cases, we indeed obtain:$$ |J_{\Phi_\Gamma}(z_1,z_2,z_3)|\sim\prod_{s=1}^3|L_1(z_s)|\left|\int_{z_1}^{z_2}\norm{\int_{z_2}^{z_3}\prod_{s=1}^2\frac{|L_2(w_s)|}{|L_1(w_s)^2|}\left|\int^{w_2}_{w_1}\frac{|L_1(y)||L_3(y)|}{|L_2(y)^2|}dy\right|dw_2}dw_1\right|.$$

Note that we can place absolute values around the $w_2$ integral as it is now clearly an integral of a real, positive integrand, and so we are done. $\hfill \qed$

So therefore,  we have, on any set in our decomposition: \begin{equation}
\begin{aligned}
\norm{J_{\Phi_\Gamma}(z_1,z_2,z_3)}\sim\prod_{s=1}^3\norm{z_s-b}^{\sigma_1}\norm{\int_{z_1}^{z_2}\norm{\int_{z_2}^{z_3}\prod_{s=1}^2\norm{w_s-b}^{\sigma_2}\left|\int^{w_2}_{w_1}\norm{y-b}^{\sigma_3}dy\right |dw_2}dw_1},
\end{aligned}
\end{equation} 
for some values of $\sigma_1,\sigma_2,\sigma_3$, depending on the type of set we are on. The explicit values are given here.
\begin{center}
 \begin{tabular}{||c| c| c| c||} 
 \hline
 &$\sigma_1$&$\sigma_2$&$\sigma_3$ \\ [0.5ex] 
 \hline\hline
 $T_{00}$ & $k_0$ & $k_{00}-2k_0 $&$ k+k_0-2k_{00}$ \\ 
 \hline
 $T_{01}$ & $0$ & $k_{01} $& $-2k_{01} $\\
 \hline
 $T_{10} $&$ k_1 $&$ k_{10}-2k_1$ &$ k_1-2k_{10}$ \\
 \hline
 $T_{11}$ &$ 0 $& $k_{11} $&$ -2k_{11}$ \\
 \hline
\end{tabular}
\end{center}
\section{The Geometric Inequality}
In this section, we establish an important inequality on the sets in our decomposition, which shall be used in both of our main theorems. This lemma is:

\begin{lemma}
For $\norm{J_{\Phi_\Gamma}}$ yielding $\sigma=(\sigma_1,\sigma_2,\sigma_3)$ values as above after the described decomposition
\begin{equation}
\begin{aligned}
\norm{J_{\Phi_\Gamma}(z_1,z_2,z_3)}\gtrsim \prod_{i=1}^3| L_\Gamma(z_i)|^\frac 13\prod_{1\leq i<j\leq 3}|z_j-z_i|,
\end{aligned}
\end{equation}

for $\sigma_3\neq -1$, and $\sigma_2+\frac{\sigma_3}{2}<-2$ or $\geq0$.

\end{lemma}
Note as this lemma is our goal, this justifies the assumption that we made in the previous section, that $c(Q)=c_1(Q)=c_2(Q)=1$. This follows, as from explicit calculation on any set, irrespective of Type, we see that, using the fact that $\norm{L_i(z)}\sim C(L_i)\norm{z-b(i)}^{\sigma_i}$, for $b(i)$ the appropriate center on this set, we have equals ``amounts'' of each of $C(L_1),$ $C(L_2)$, and $C(L_3)$ on either side of (3.1), so they cancel out, and as such, establishing this lemma under our assumptions is indeed satisfactory to prove the general case.

Furthermore, we can meet the conditions of this lemma, that is to say we can avoid certain values of $\sigma$, by applying an affine transformation to our curve $\Gamma$ before applying the decomposition $D2$. This is done to ensures that we avoid certain $\sigma$ values, as described in \cite{DendrinosWright}. Specifically, we ensure there is no $k_{10}$ value such that $k_{10}=\frac{-1-k_1}2$ on $T_{10}$ intervals, and no $k_{00}$ value such that $k_{00}=\frac{k+k_0+1}2$ on $T_{00}$ intervals. This ensures $\sigma_3\neq-1$. Similarly, we avoid $k_1=1$ on $T_{10}$ and  and any $k_0$ values on $T_{00}$ such that $k_0=\frac{k+i}{3}$ for $i=1,2,3,4$. Note that as $k_0$ is an integer, this involves excluding at most $2$ values of $k_0$.

\emph{Proof of Lemma 3.1:}
To establish the lemma, we first seek to show that: 
\begin{equation}
\begin{aligned}
\norm{J_{\Phi_\Gamma}(z_1,z_2,z_3)} \gtrsim \prod_{i=1}^3\norm{z_i}^{\sigma_1+\frac{2\sigma_2}{3}+\frac{\sigma_3}{3}}\norm{\int_{z_1}^{z_2}\int_{z_2}^{z_3}\norm{\int_{w_1}^{w_2}dy}dw_1dw_2}.
\end{aligned}
\end{equation}

Again, calculation will yield that $\norm{z_i}^{\sigma_1+\frac{2\sigma_2}{3}+\frac{\sigma_3}{3}}\sim\norm{L_\Gamma(z_i)}^{\frac 13}$. 

To establish this, we first wish to see that: $$I_1:=\norm{\int^{w_2}_{w_1}\norm{y-b}^{\sigma_3}dy}\gtrsim\norm{w_2-b}^{\frac{\sigma_3}{2}}\norm{w_1-b}^{\frac{\sigma_3}{2}}\norm{\int_{w_1}^{w_2}dy}$$

Firstly, we may assume $\norm{w_2}>\norm{w_1}$, as if not, we can notice that

$\norm{\int^{w_2}_{w_1}\norm{y-b}^{\sigma_3}dy}=\norm{-\int^{w_1}_{w_2}\norm{y-b}^{\sigma_3}dy}=\norm{\int^{w_1}_{w_2}\norm{y-b}^{\sigma_3}dy}$.

So we may relabel to ensure that the larger of the two occurs at the upper bound of the integral.

Also notice that if $\norm{w_1-b}<\norm{w_2-b}<9\norm{w_1-b}$, we have  $\norm{w_1-b}\sim\norm{y-b}\sim\norm{w_2-b}$ throughout this line integral, and our result immediately follows. We therefore only have to deal with the case in which $9\norm{w_1-b}<\norm{w_2-b}$. 

Let: \begin{itemize}
    \item $L:=\{y\mid2\norm{w_1-b}<\norm{y-b}<3\norm{w_1-b}\}.$ 
    \item $U:=\{y \mid\frac13\norm{w_2-b}<\norm{y-b}<\frac12\norm{w_2-b}\}.$
\end{itemize}. Note that after parameterizing our line integral, we can see that it is the integral of a purely positive quantity over a real interval, so we can delete portions of the line segment and ensure that this results in a lower bound.

So 
\begin{equation*} 
\begin{split}
I_1&\geq\norm{\left[\int_L+\int_U\right]{\norm{y-b}^{\sigma_3}}dy}\\
&\sim\norm{w_1-b}^{\sigma_3+1}+\norm{w_2-b}^{\sigma_3+1}\\
&=\norm{w_2-b}\left(\norm{w_2-b}^{\sigma_3}+\frac{\norm{w_1-b}^{\sigma_3+1}}{\norm{w_2-b}}\right).
\end{split}
\end{equation*}

Note that $\norm{w_2-b}\sim\norm{w_2-w_1}=\norm{\int_{w_1}^{w_2}dy}$, as $9\norm{w_1-b}<\norm{w_2-b}$

So $ I_1\gtrsim\norm{\int_{w_1}^{w_2}dy}\text{max}\left(\norm{w_2-b}^{\sigma_3},\frac{\norm{w_1-b}^{\sigma_3+1}}{\norm{w_2-b}}\right) $

If $\sigma_3$ is positive, we have $$I_1\gtrsim\norm{\int_{w_1}^{w_2}dy}\norm{w_2-b}^{\sigma_3}>\norm{\int_{w_1}^{w_2}dy}\norm{w_2-b}^{\frac{\sigma_3}{2}}\norm{w_1-b}^{\frac{\sigma_3}{2}}.$$

If $\sigma_3$ is negative, then, as $\sigma_3$ is an integer which is not equal to -1, $\sigma_3+1$ is negative also, so $$I_1\gtrsim\norm{\int_{w_1}^{w_2}dy}\frac{\norm{w_1-b}^{\sigma_3+1}}{\norm{w_2-b}}>\norm{\int_{w_1}^{w_2}dy}\norm{w_2-b}^{\frac{\sigma_3}{2}}\norm{w_1-b}^{\frac{\sigma_3}{2}}.$$

As this holds for all $\sigma_3\neq-1$, we are done. With this, we can conclude that$$
\norm{J_{\Phi_\Gamma}(z_1,z_2,z_3)}\gtrsim\prod_{s=1}^3\norm{z_s-b}^{\sigma_1}\norm{\int_{z_1}^{z_2}\int_{z_2}^{z_3}\prod_{s=1}^2\norm{w_s-b}^{\sigma_2+\frac{\sigma_3}{2}}\norm{\int_{w_1}^{w_2}dy}dw_2dw_1}.
$$

We now seek to establish a similar inequality for $$I_2=\norm{\int_{z_2}^{z_{3}}\norm{w_2-b}^{\sigma_2+\frac{\sigma_3}{2}}\norm{\int_{w_1}^{w_2}dy}dw_2}.$$

We want to see that $I_2\gtrsim\text{max}\left(\norm{z_{3}-b}^{\sigma_2+\frac{\sigma_3}{2}},\frac{\norm{z_{2}-b}^{\sigma_2+\frac{\sigma_3}{2}+2}}{\norm{z_{3}-b}^{2}}\right)\norm{\int_{z_2}^{z_{3}}\norm{\int_{w_1}^{w_2}dy}dw_2}$,\newline for $\sigma_2+\frac{\sigma_3}{2}<-2$ or $\geq 0$.

We assume that $\norm{z_2-b}<\norm{z_{3}-b}$, and relabel if we are not in such a case. We also notice that the inequality is immediate if we are in the case that $\norm{z_2-b}<\norm{z_{3}-b}<9\norm{z_2-b}$, so we therefore assume that $9\norm{z_2-b}<\norm{z_{3}-b}$, and again restrict the integral to $L$ and $U$, which now represent the sets:
\begin{itemize}
    \item $L:=\{y\mid2\norm{z_2-b}<\norm{y-b}<3\norm{z_2-b}\}$.
    \item $U:=\{y \mid\frac13\norm{z_{3}-b}<\norm{y-b}<\frac12\norm{z_{3}-b}\}$.
\end{itemize}
So  we have $I_2\geq\norm{\left[\int_L+\int_U\right]\norm{w_2-b}^{\sigma_2+\frac{\sigma_3}{2}}\norm{w_2-w_1}dw_2}$. Identically to last time, we have $\norm{w_2-w_1}\sim\norm{w_2-b}$, so we can write:


\begin{equation*}
\begin{split}
I_2 & \gtrsim\norm{\left[\int_L+\int_U\right]\norm{w_2-b}^{\sigma_2+\frac{\sigma_3}{2}+1}dw_2}\\
 & \sim \norm{z_2-b}^{\sigma_2+\frac{\sigma_3}{2}+2}+\norm{z_3-b}^{\sigma_2+\frac{\sigma_3}{2}+2}\\
 &=\norm{z_3-b}^2\left(\frac{\norm{z_2-b}^{\sigma_2+\frac{\sigma_3}{2}+2}}{\norm{z_3-b}^2}+\norm{z_3-b}^{\sigma_2+\frac{\sigma_3}{2}}\right)\\
 &\gtrsim\text{max}\left(\norm{z_{3}-b}^{\sigma_2+\frac{\sigma_3}{2}},\frac{\norm{z_{2}-b}^{\sigma_2+\frac{\sigma_3}{2}+2}}{\norm{z_{3}-b}^{2}}\right)\norm{\int_{z_2}^{z_{3}}\norm{\int_{w_1}^{w_2}dy}dw_2}.
\end{split}
\end{equation*}

Now, if $\sigma_2+\frac{\sigma_3}{2}>0$, we have:

\begin{equation*}
\begin{split}
I_2 & = \norm{z_{3}-b}^{\sigma_2+\frac{\sigma_3}{2}}\norm{\int_{z_2}^{z_{3}}\norm{\int_{w_1}^{w_2}dy}dw_2}\\
 & \geq\norm{z_{3}-b}^{\frac{2\sigma_2}{3}+\frac{\sigma_3}{3}}\norm{z_{2}-b}^{\frac{\sigma_2}{3}+\frac{\sigma_3}{6}}\norm{\int_{z_2}^{z_{3}}\norm{\int_{w_1}^{w_2}dy}dw_2}.
\end{split}
\end{equation*}

if instead, $\sigma_2+\frac{\sigma_3}{2}<-2$, we have

\begin{equation*}
\begin{split}
I_2&\gtrsim\frac{\norm{z_{2}-b}^{\sigma_2+\frac{\sigma_3}{2}+2}}{\norm{z_{3}-b}^{2}}\norm{\int_{z_2}^{z_{3}}\norm{w_2-w_1}dw_2}\\
 & >\norm{z_{3}-b}^{\frac{2\sigma_2}{3}+\frac{\sigma_3}{3}}\norm{z_{2}-b}^{\frac{\sigma_2}{3}+\frac{\sigma_3}{6}}\norm{\int_{z_2}^{z_{3}}\norm{\int_{w_1}^{w_2}dy}dw_2},
\end{split}
\end{equation*}

as  $\sigma_2+\frac{\sigma_3}{2}+2<0$
 
A similar analysis yields that for $$I_3=\norm{\int_{z_1}^{z_{2}}\norm{w_1-b}^{\sigma_2+\frac{\sigma_3}{2}}\norm{\int_{w_1}^{w_2}dy}dw_1},$$
 
 we have $$ I_3\gtrsim\norm{z_{1}-b}^{\frac{2\sigma_2}{3}+\frac{\sigma_3}{3}}\norm{z_{2}-b}^{\frac{\sigma_2}{3}+\frac{\sigma_3}{6}}\norm{\int_{z_2}^{z_{3}}\norm{\int_{w_1}^{w_2}dy}dw_2},$$ for the same constraints on our $\sigma$ values as previously outlined. So we therefore indeed have (3.2). 
 
 From here we seek to resolve the triple integral which appears in (3.2), and establish that:
 \begin{equation}
\begin{aligned}
\norm{\int_{z_1}^{z_2}\norm{\int_{z_2}^{z_3}\left|\int_{w_1}^{w_2}dy\right|dw_2}dw_1}\gtrsim \norm{z_3-z_1}\norm{z_3-z_2}\norm{z_2-z_1}.
\end{aligned}
\end{equation}
Clearly the inner most integral is given by $\norm{\int^{w_2}_{w_1}dy}=\norm{w_2-w_1}$.

Consider  $J=\norm{\int_{z_2}^{z_{3}}\norm{w_2-w_1}dw_2}$. If $z_3-z_2=\norm{z_3-z_2}e^{i\phi}$, and $t=e^{-i\phi}(w_2-z_2)\in\R$, we have $$J=\int_0^{\norm{z_3-z_2}}\norm{t-w_1'}=\int_0^{\norm{z_3-z_2}}\sqrt{(t-\text{Re}(w_1'))^2+\text{Im}(w_1')^2}dt:=\int_0^{\norm{z_3-z_2}}g_1(t)dt,$$
 where $w_1'=e^{-i\phi}(w_1-z_2)$.
 
 To bound this from below, we introduce a new function $f_1(t)=\norm{t-\widehat{w_1'}(t)}$, where    $$\widehat{w_1'}(t)= 
\begin{dcases}
    w_1'-i\frac{\text{Im}(w_1')}{\text{Re}(w_1')}t& \text{if } t\leq \text{Re}(w_1')\\
    w_1'-i\text{Im}(w_1')\frac{\norm{z_3-z_2}-t}{\norm{z_3-z_2}-\text{Re}(w_1')},              & t>\text{Re}(w_1')
\end{dcases}.
$$
Note that, for all values of $t\in[0,|z_3-z_2|]$, we have that $\text{Re}(w_1')=\text{Re}(\widehat{w_1'}(t))$, and also that $|\text{Im}(w_1')|>|\text{Im}(\widehat{w_1'}(t)|$. We also have $\widehat{w_1'}(0)=\widehat{w_1'}(|z_3-z_2|)=w_1'$.

Clearly $g_1(t)\geq f_1(t)$ for all $t$ in the interval $[0,|z_3-z_2|]$, as $t$ is real, and $\widehat{w_1'}(t)$ is  $w_1'$ with a smaller imaginary part. Note that computation shows that $f_1(t)$ is given by :

$$f_1(t)= 
\begin{dcases}
    \frac{\norm{z_2-w_1}}{\text{Re}(w_1')}(\text{Re}(w_1')-t)& \text{if } t\leq \text{Re}(w_1')\\
    \frac{\norm{z_3-w_1}}{\norm{z_3-z_2}-\text{Re}(w_1')}(t-\text{Re}(w_1')) & t>\text{Re}(w_1')
\end{dcases}.
$$

So $J=\int_0^{\norm{z_3-z_2}}g_1(t)dt>\int_0^{\norm{z_3-z_2}}f_1(t)dt$. Now note that $g_1(t)$ increases as $\norm{t-\text{Re}(w_1')}$ increases. Also $g_1(0)=\norm{z_2-w_1}=f_1(0)$ and $g_1(\norm{z_3-z_2})=\norm{z_3-w_1}=f_1(\norm{z_3-z_2})$. 

Therefore if $\text{Re}(w_1')<\frac{\norm{z_3-z_2}}{2}$, we have $\norm{z_3-w_1}>\norm{z_2-w_1}$.

Alternatively if $\text{Re}(w_1')>\frac{\norm{z_3-z_2}}{2}$, then we have $\norm{z_2-w_1}>\norm{z_3-w_1}$.

Now, suppose $\text{Re}(w_1')\leq0$.
Then $\norm{z_3-w_1}>\norm{z_2-w_1}$, and we have:

\begin{equation*}
\begin{split}
J&>\int_0^{\norm{z_3-z_2}}f_1(t)dt\\
&=\int_0^{\norm{z_3-z_2}} \frac{\norm{z_3-w_1}}{\norm{z_3-z_2}-\text{Re}(w_1')}(t-\text{Re}(w_1'))dt\\
 & = \frac{\norm{z_3-w_1}}{\norm{z_3-z_2}-\text{Re}(w_1')}\left[\frac{t(t-\text{Re}(w_1'))}{2}\right]^{\norm{z_3-z_2}}_0\\
 &\sim\norm{z_3-z_2}\norm{z_3-w_1}.
\end{split}
\end{equation*}

A symmetric argument gives that if $\text{Re}(w_1')>\norm{z_3-z_2}$
then $\norm{z_2-w_1}>\norm{z_3-w_1}$, and $J\gtrsim\norm{z_3-z_2}\norm{z_2-w_1}$.

Lastly, if $\text{Re}(w_1')\in\left[0,\norm{z_3-z_2}\right]$, then:
\begin{equation*} 
\begin{split}
J&>\int_0^{\text{Re}(w_1')}\frac{\norm{z_2-w_1}}{\text{Re}(w_1')}(\text{Re}(w_1')-t)dt+\int_{\text{Re}(w_1')}^{\norm{z_3-z_2}}\frac{\norm{z_3-w_1}}{\norm{z_3-z_2}-\text{Re}(w_1')}(t-\text{Re}(w_1'))dt\\
 & =\frac12\left(\norm{z_2-w_1}\text{Re}(w_1')+\norm{z_3-w_1}(\norm{z_3-z_2}-\text{Re}(w_1') )\right)\\
 &\sim\norm{z_3-z_2}\text{max}\left(\norm{z_3-w_1},\norm{z_2-w_1}\right).
\end{split}
\end{equation*}

So indeed $J\gtrsim\sim\norm{z_3-z_2}\text{max}(\norm{z_3-w_1},\norm{z_2-w_1})$ in each of these cases. Applying this to the initial double integral, we get $$\norm{\int_{z_1}^{z_2}\int_{z_2}^{z_{3}}\norm{w_2-w_1}dw_2dw_1}\gtrsim\norm{z_3-z_2}\norm{\int_{z_1}^{z_2}\norm{z_3-w_1}dw_1} $$
$$\gtrsim\norm{z_3-z_2}\norm{z_2-z_1}\norm{z_3-z_1}.$$

As we can go through an identical procedure for the remaining single integral. 

Therefore, we have established (3.3). This in conjunction with (3.2) yields the proof of (3.1). \qed
\section{Application to the Convolution Operator}
We will now use the work from Sections 2 and 3 to establish Theorem 1.1. Let us denote the measure $d\sigma(z):=\lambda_\gamma(z)dz.$

In order to prove Theorem 1.1, and Corollary 1.2,  we need only establish that the operator $T$ is of restricted weak-type at $(p,q)=(2,3)$. The bound at $(0,0)$ follows trivially from the fact that $T_\Gamma f$ is bounded whenever $f$ is bounded. Furthermore, as the operator is ``almost'' self dual, we can use duality properties to conclude that we also have the bound corresponding to the vertex $(1,1)$.

Therefore we only need to establish the restricted weak-type bound corresponding to the vertex $\left(\frac 12,\frac 13\right)$, as this will give us the bound at $\left(\frac 23,\frac 12\right)$, by duality. From here, we can use the Marcinkiewicz Interpolation Theorem between these four vertices to establish our desired range, the convex hull of these points, $\mathcal{R}$. The bound we establish at $\left(\frac 12,\frac 13\right)$ is uniform, as opposed to the trivial bound at $(0,0)$. Therefore, we have uniform estimates along the line segment joining $\left(\frac 12,\frac 13\right)$ and $\left(\frac 23,\frac 12\right)$, and local estimates, with implicit constants dependent on the radius, at every other $(p,q)$ pair in $\mathcal{R}$.

To establish the restricted weak-type, we introduce two quantities, $\alpha$ and $\beta$, where, for measurable sets $E$ and $F$, they are given by the equations: \begin{equation}
\begin{aligned}
\alpha\norm{F}=\langle T{_\Gamma} \chi_{_E}, \chi_{_F}\rangle=\beta\norm{E}.
\end{aligned}
\end{equation}

For $\chi_S$ being the characteristic function of the set $S$. Our estimate at $(2,3)$ is equivalent to the statement that for all measurable sets, $E$ and $F$, we have: 
$$
\langle T{_\Gamma} \chi_{_E}, \chi_{_F}\rangle\lesssim \norm{E}^\frac{1}{2}\norm{F}^\frac{2}{3}.$$

which can be rewritten as
\begin{equation}
\begin{aligned}
\norm{E}\gtrsim\alpha^4\beta^2.
\end{aligned}
\end{equation}

which is the statement we will prove. To do this, we construct subsets of $E$ and $F$ as follows.

Let $E_1:=\{y\in E\:|\: T^*\chi_F(y)\geq\frac{\beta}{2}\}$. Note that: $$\langle T\chi_{E_1},\chi_F \rangle=\langle T\chi_{E},\chi_F \rangle-\langle T\chi_{E\setminus E_1},\chi_F \rangle\geq\alpha\norm{F}-\frac{\beta}{2}\norm{E}=\frac\alpha 2\norm{F}.$$

That is to say that the average value of $T\chi_{E_1}$ on $F$ is greater than $\frac\alpha 2$. We therefore can introduce the set $F_1:=\{x\in F\:|\:T\chi_{E_1}>\frac\alpha 2\}$, and this will be non-empty.
Take $z_0\in F_1$, and let $P:=\{z_1\:|\:z_0-\Gamma(z_1)\in E_1\}$. Then we have $\sigma(P)=T\chi_{E_1}(z_0)>\frac\alpha 4$.

Let $Q_{z_1}=\{z_2\:|\:z_0-\Gamma(z_1)+\Gamma(z_2)\in F_1\}$. Again, $\sigma(Q_{z_1})=T^*\chi_F(z_0-\Gamma(z_1))\geq\frac\beta 2$

Let $R_{z_1,z_2}=\{z_3\:|\:z_0-\Gamma(z_1)+\Gamma(z_2)-\Gamma(z_3)\in E_1\}$, and identically to the other sets, we have that $\sigma(R_{z_1,z_2})=T\chi_{E_1}(z_0-\Gamma(z_1)+\Gamma(z_2))\geq\frac\alpha 2$. 

We let $S=P\times Q_{z_1}\times R_{z_2,z_2}$. 

We have $z_0+\Phi_\Gamma(S) \subseteq E_1\subseteq E$, where $\Phi_\Gamma(z_1,z_2,z_3)=-\Gamma(z_1)+\Gamma(z_2)-\Gamma(z_3)$. Furthermore, we restrict ourselves to the subset of $S$ on which $\norm{z_1}\leq\norm{z_2}\leq\norm{z_3}$ to ensure $\Phi_\Gamma$ is c(N) to 1 on the region, and we may apply Bézout's Theorem. By Bézout's Theorem, we have, for $J_\R$, and $J_\C$ being the real and complex jacobian of the map of $\Phi_\Gamma$, the following estimates:

\begin{equation} 
\begin{split}
\norm{E} & \gtrsim \int_S \norm{J_\R(z_1,z_2,z_3)}dz_3dz_2dz_1 = \int_S \norm{J_\C(z_1,z_2,z_3)}^2dz_3dz_2dz_1 \\
 & \gtrsim\int_S\prod_{i=1}^3|z_i|^{2\sigma_1+\frac{4\sigma_2}{3}+\frac{2\sigma_3}{3}}\prod_{i=1,j<i}^3|z_i-z_j|^2dz_3dz_2dz_1.
\end{split}
\end{equation}

Therefore, for $k'=k$ or $=0$ appropriately, we have: $$\norm{E}\gtrsim \int_S\ \norm{z_1}^{\frac{2k'}{3}}\norm{z_2}^{\frac{2k'}{3}}\norm{z_3}^{\frac{2k'}{3}}\prod_{i=1,j<i}^3\norm{z_i-z_j}^2dz_3dz_2dz_1.$$
First consider the set $B_x=\{z\in \C\:|\:\norm{z}<(16\pi\nu)^{-\nu}x^\nu\}$, for $\nu=\frac{3}{k'+6}$. \newline
Now $\sigma(B_x)=\int_{B_x}\norm{z}^\frac {k'}{3}dz=2\pi\int_0^{(16\pi\nu)^{-\nu}x^\nu}r^{\frac {k'}{3} +1}dr=\frac{2\pi}{\frac {k'}{3} +2}((16\pi\nu)^{-\nu}x^\nu)^{\frac {K'}{3}+2}=\frac x8$

So therefore, we can delete $B_\alpha$ from $P$ and $R_{z_1,z_2}$, and $B_\beta$ from $Q_{z_1}$ and keep these sets having sigma-measure comparable to $\alpha$ or $\beta$ as appropriate.

Consider the case where two of these variables are of comparable modulus. Let us suppose we have $\norm{z_j}<\norm{z_i}<2\norm{z_j}$, for $j<i$. We also assume, without loss of generality, that $z_i$ is in a set with sigma-measure comparable to $\alpha$. That is to say, that $i$ is odd.

Now suppose $z_j\in B_\frac\alpha 2$. In this case, we have $\norm{z_i-z_j}>\norm{z_i}-\norm{z_j}>\norm{z_i}(1-\frac{1}{2^\nu})$. Furthermore, we have $ \norm{z_i}>(16\pi\nu)^{-\nu}\alpha^\nu$. It follows $\norm{z_i}^\frac{k'+6}{6}>(16\pi\nu)^{-\frac{1}{2}}\alpha^\frac 12$, which gives us $\norm{z_i}\gtrsim\alpha^\frac 12\norm{z_i}^{-\frac {k'}{6}}$. So therefore, $\norm{z_i-z_j}\gtrsim\beta^\frac 12(\norm{z_i}\norm{z_j})^{-\frac{k'}{12}}$, as $\norm{z_i}\sim\norm{z_j}$.

In the case where $z_j\notin B_\frac {\alpha}{2}$, we have $\norm{z_j}>(16\pi\nu)^{-\nu}\left(\frac{\alpha}{2}\right)^\nu$. It follows, similar to the other case, that $\norm{z_j}> (32\pi\nu)^{-\frac 12}\alpha^\frac 12\norm{z_j}^{-\frac{k'}{6}}$. For $S_i=P$ or $=R_{z_1,z_2}$ depending on our $i$ value, consider $B_\alpha(z_j)=\{z_i\in S_i|\norm{z_i-z_j}<c_0\norm{z_j}^{-\frac {k'}{6}}\alpha^\frac 12\}$.

Now on $B_\alpha(z_j)$, we have $\norm{z_i-z_j}<c_0\norm{z_j}^{-\frac {k'}{6}}\alpha^\frac 12<c_0(32\pi\nu)^\frac 12\norm{z_j}$. Consider the sigma-measure of this set. We have $\sigma(B_\alpha(z_j))=\int_{B_\alpha(z_j)}\norm{z}^\frac {k'}{3} dz$. But for all $z\in B_\alpha(z_j)$, we have:\begin{equation*}
\norm{z}\leq\norm{z_j}+\norm{z-z_j}\leq\norm{z_j}(1+c_0(32\pi\nu)^\frac 12).
\end{equation*}
Therefore, $\sigma(B_\alpha(z_j))<\norm{z_j}^\frac {k'}{3}(1+c_0(32\pi\nu)^\frac 12)^\frac {k'}{3}\pi(c_0\norm{z_j}^{-\frac {k'}{6}}\alpha^\frac 12)^2$. The right hand side of this inequality is comparable to $\alpha$, for $c_0$ chosen sufficiently small. It can therefore be removed from $S_i$, and still have $\sigma(S_i)\sim \alpha$. Deleting this set gives us $\norm{z_i-z_j}\gtrsim \alpha^\frac 12 (\norm{z_i}\norm{z_j})^{-\frac {k'}{12}}$, as in the other case.

Therefore, for $i$ odd, $j<i$, if $\norm{z_i}\sim\norm{z_j}$, we have $\norm{z_i-z_j}\gtrsim\alpha^\frac12(\norm{z_i}\norm{z_j})^{-\frac {k'}{12}}$. An identical argument for the $i=2$  case gives us $\norm{z_i-z_j}\gtrsim\beta^\frac12(\norm{z_i}\norm{z_j})^{-\frac {k'}{12}}$.

Also notice that if we are not in this case, that is we do not have $\norm{z_i}\sim \norm{z_j}$, we necessarily have $\norm{z_i}>2\norm{z_j}$, which gives $\norm{z_i-z_j}\sim \norm{z_i}$. Using these facts, our ordering of the size of $z_1$, $z_2$, and $z_3$,  and our integral bound for $\norm{E}$, we can establish $\norm{E}\gtrsim\alpha^4\beta^2$ as desired. We shall prove as an example the case in which $\norm{z_3}\sim\norm{z_2}>2\norm{z_1}$, and the other cases can be dealt with similarly.
\begin{align*}
\norm{E}&\gtrsim \int_S\ \norm{z_1}^{\frac{2k'}{3}}\norm{z_2}^{\frac{2k'}{3}}\norm{z_3}^{\frac{2k'}{3}}\prod_{i=1,j<i}^3\norm{z_i-z_j}^2dz_3dz_2dz_1\\
&=\int_S\ \norm{z_1}^{\frac{k'}{3}}\norm{z_2}^{\frac{k'}{3}}\norm{z_3}^{\frac{k'}{3}}\prod_{i=1,j<i}^3\norm{z_i-z_j}^2\norm{z_1}^{\frac{k'}{3}}\norm{z_2}^{\frac{k'}{3}}\norm{z_3}^{\frac{k'}{3}}dz_3dz_2dz_1\\
&\sim\int_S\ \norm{z_1}^{\frac{k'}{3}}\norm{z_2}^{\frac{k'}{3}}\norm{z_3}^{\frac{k'}{3}}\
 \norm{z_3-z_2}^2\norm{z_3}^2\norm{z_2}^2
 \norm{z_1}^{\frac{k'}{3}}\norm{z_2}^{\frac{k'}{3}}\norm{z_3}^{\frac{k'}{3}}dz_3dz_2dz_1\\
 &\gtrsim\int_S\ \norm{z_1}^{\frac{k'}{3}}\norm{z_2}^{\frac{k'}{3}}\norm{z_3}^{\frac{k'}{3}}\
(\alpha^\frac12(\norm{z_3}\norm{z_2})^{-\frac {k}{12}})^2
 \norm{z_1}^{\frac{k'}{3}}\norm{z_2}^{\frac{k'}{3}+2}\norm{z_3}^{\frac{k'}{3}+2}dz_3dz_2dz_1\\
 &\gtrsim\int_S\ \norm{z_1}^{\frac{k'}{3}}\norm{z_2}^{\frac{k'}{3}}\norm{z_3}^{\frac{k'}{3}}\
(\alpha^\frac12\norm{z_3}^{-\frac {k}{6}})^2
 \norm{z_1}^{\frac{k'}{3}+2}\norm{z_2}^{\frac{k'}{3}+2}\norm{z_3}^{\frac{k'}{3}}dz_3dz_2dz_1\\
 &=\alpha\int_S\ \norm{z_1}^{\frac{k'}{3}}\norm{z_2}^{\frac{k'}{3}}\norm{z_3}^{\frac{k'}{3}}\
 \norm{z_1}^{\frac{k'}{3}+2}\norm{z_2}^{\frac{k'}{3}+2}dz_3dz_2dz_1\\
 &=\alpha\int_S\ \norm{z_1}^{\frac{k'}{3}}\norm{z_2}^{\frac{k'}{3}}\norm{z_3}^{\frac{k'}{3}}\
((16\pi\nu)^{-2\nu}(\alpha^\nu\beta^\nu))^{\frac{k'}{3}+2}dz_3dz_2dz_1\\
&\sim\alpha^2\beta\int_P\int_{Q_{z_1}}\int_{R_{z_1.z_2}}\ \norm{z_1}^{\frac{k'}{3}}\norm{z_2}^{\frac{k'}{3}}\norm{z_3}^{\frac{k'}{3}}dz_3dz_2dz_1\\
&\sim\alpha^2\beta\sigma(P)\sigma(Q_{z_1})\sigma(R_{z_1,z_2})\\
&\sim\alpha^4\beta^2.
\end{align*}
This can be achieved in all cases. The general strategy is to first split the monomials into a product of two parts of equal size, that is writing them as  $\norm{z_i}^\frac{2k'}{3}=\norm{z_i}^\frac{k'}{3}\norm{z_i}^\frac{k'}{3}$. We aim to keep one of these in the integrand, in order to eventually turn our integral into the sigma measure of our sets, and the other will be used to generate either an $\alpha$ or a $\beta$.

Then, for any comparable-size variables, $z_i$ and $z_j$, $i>j$, we apply our bounds for $\norm{z_i-z_j}$, this generates one $\alpha$ or $\beta$, and a $\norm{z_i}^{-\frac{k'}{3}}$, which we use to remove one of the two parts we split our initial monomial of that variable into. 

This may still leave us with some of the monomials that we still wish to replace with $\alpha$ or $\beta$.
To do this, for any variables that are not comparable, $z_i$ and $z_j$, $i>j$, we replace $\norm{z_i-z_j}^2$ with $\norm{z_i}^2$. From here we can redistribute the exponents of the remaining variables ``downwards'', that is to variables with lower subscripts, so that all our remaining monomials are now monomials of the form $\norm{z_i}^{\frac{k'}{3}+2}$. We can then use our lower bounds that we get from deleting $B_\alpha$ or $B_\beta$ to replace this with quantities comparable to $\alpha$ or $\beta$ as appropriate. $\hfill\qed$

Note that as the estimate here is uniform in $r$, we can in fact extend the restricted weak-type estimate to the entirety of $\C$ by a simple re-scaling argument and the Monotone Convergence Theorem. We therefore have global $(p,q)$ estimates for $T_\Gamma$ for all $(p^{-1},q^{-1})$ on the line segment joining $\left(\frac 12,\frac 13\right)$, and $\left(\frac 23,\frac 12\right)$.
\newpage
\section{Application to Fourier Restriction}

For Theorem 1.3, we again rely on interpolation to achieve the desired range. Again, it is trivial to establish the estimate corresponding to vertex $(1,0)$, so we must only establish the bound at $(\frac67,\frac 67)$. To achieve this, we will use the dual estimate. 
The dual to estimate in this theorem is given by, for $q>7$, and $p'=\frac q6$, where $p'$ is the Hölder Conjugate of $p$: \begin{equation}
\begin{aligned}
||\mathcal{E}_\Gamma f||_{L^q}\leq C||f||_{L^p(\lambda_\Gamma(z)dz}.
\end{aligned}
\end{equation}

Here, we define the Extension Operator associated to $\Gamma$ as: 
\begin{equation}
\begin{aligned}
\mathcal{E}_\Gamma f(\mathbf{z})=\int_\C e^{i\mathbf{z}\cdot\Gamma(w)}f(w)\lambda_\Gamma(w)dw.
\end{aligned}
\end{equation}

Now to establish this, we will use the method developed by Stovall in \cite{Stovall}, which establishes the analogous result in the real case. This is more involved that the convolution estimate, but mostly follows Stovall's real method closely. As in that paper, we first establish a bound in the case where we are on a set on which our torsion $L_\Gamma(z)$ is bounded. This appears as Theorem 2.1 in \cite{Stovall}. It translates into the complex case as follows:

\begin{lemma}
Suppose we have $(p,q)$ satisfying $q>7$, and $p'=6q$, a degree $N$ curve $\gamma(z)$, and a convex set $B\subset\C$ such that for all $z\in B$, we have that $$0<C_1<|L_\gamma(z)|<C_2,$$ then 
$$||\mathcal{E}_\gamma\chi_B f||_{L^q}\lesssim||f||_{L^p(\lambda_\gamma(z)dz)},$$
with implicit constants depending on $N,p$ and the ratio $\frac {C_2}{C_1}$
\end{lemma}

Note, as $L_\gamma(z)$ is a polynomial of degree at most $3N$, that we can decompose $B=\bigcup_{j=1}^{3\log(\frac{C_2}{C_1})N}B_j$ where on each $B_j$ we have $\frac C2<|L_\gamma(z)|<2C $, and by the triangle inequality, we may assume that our initial $B$ is replaced with one of these $B_j$. 

Furthermore, as our estimate is affine invariant, we may apply an affine transformation to assume that $C=1$, and that  $B$ contains some fixed complex number, say $0$ in it's interior, and $|L_\gamma(0)|\neq 0$.

To prove this theorem, we must first prove the analogue to Lemma 2.2 in \cite{Stovall}, which translates directly into the Complex case without notable changes, and is only included for completion, and reads as:

\begin{lemma}
Suppose we have a curve $\gamma$, and a the set $B$ satisfying the previous theorems conditions. Then there exists an affine transformation $A$, such that $\det A=1$, and $||A\gamma ||_{C^N(B)}\lesssim 1$, with implicit constants depending only on $N$.
\end{lemma}

\emph{Proof of Lemma:} 
As in \cite{Stovall}, we take $A$ such that $Az=\left(\gamma'(0),\gamma''(0),\gamma'''(0)^{-1}\right)(z-\gamma(0))$. As $0\in B$, and understanding $\det A$ refers to the determinant of the matrix component of $A$, we have $\det A\in [\frac 12,2]$. Thus to prove the lemma, replacing $A$ with $\frac 1 {\det A} A$ if necessary, we must just show $||A\gamma ||_{C^N(B)}\lesssim 1$.

Note that after this transformation, we have the following properties for $A\gamma$:

\begin{itemize}
    \item $A\gamma(0)=0$.
    \item $|L_{A\gamma}(z)|\in [\frac 14,4]$ on $B$.
    \item For $j\in \{1,2,3\}$, we have $A\gamma^{(j)}(0)=e_j$.
\end{itemize}

From here, we will use the notation that $\gamma=A\gamma$. Note that by Taylor's Theorem, we need only estimate the size of the remaining $\gamma^{(j)}(0)$ for $j>3$, or, sufficiently, by the properties of $\gamma$, that for all $(n_1,n_2,n_3)$, with $n_1<n_2<n_3<N$, we have $$\text{det}\left(\gamma^{(n_1)}(0),\gamma^{(n_2)}(0),\gamma^{(n_3)}(0)\right)\lesssim 1.$$

To see this, we assume a contradiction. That is to say there is a sequence $\{\gamma_n\}_{n=1}^\infty$ with the properties listed above, and that $$\lim_{n\rightarrow\infty} \max_{n_1<n_2<n_3<N}\det\left(\gamma_n^{(n_1)}(0),\gamma_n^{(n_2)}(0),\gamma_n^{(n_3)}(0)\right) =\infty.$$

Note that for the maximising $(n_1,n_2,n_3)$, we have $n_1+n_2+n_3>6$, for all sufficiently large $n$, as otherwise the above determinant is $1$.

Now introduce a sequence of real numbers $\{\delta_n\}_{n=1}^\infty$ with $\delta_n\in(0,1)$, and for each $\gamma_n$, let $\Gamma_n$ be the curve such that its components are given by $\Gamma_{n,j}(z)=\delta_n^{-j}\gamma_{n,j}(\delta_n z)$.

Note that $\Gamma_n$ satisfies the properties listed above on the set $B_{\delta_n}$, which is the set $B$ scaled from $0$ by a factor of $\frac 1 {\delta}$ Furthermore, calculation yields for $(n_1,n_2,n_3)\neq(1,2,3)$, $$\det\left(\Gamma_n^{(n_1)}(0),\Gamma_n^{(n_2)}(0),\Gamma_n^{(n_3)}(0)\right)=\delta_n \det\left(\gamma_n^{(n_1)}(0),\gamma_n^{(n_2)}(0),\gamma_n^{(n_3)}(0)\right)$$
Therefore, for sufficiently large $n$, we can choose $\delta_n$ such that $$\max_{n_1<n_2<n_3<N}\det\left(\Gamma_n^{(n_1)}(0),\Gamma_n^{(n_2)}(0),\Gamma_n^{(n_3)}(0)\right) =5.$$

Note that this implies that $\{\delta_n\}_{n=1}^\infty$ converges to $0$. We can also take a subsequence of $\Gamma_n$ such that for this  new sequence: $$\det\left(\Gamma_n^{(n_1)}(0),\Gamma_n^{(n_2)}(0),\Gamma_n^{(n_3)}(0)\right) =5,$$

for all $n$, and for some fixed $(n_1,n_2,n_3)$. This can be done as there are only finitely many choices for $n_1,n_2,$ and $n_3$.

However, this paired with the properties of $\Gamma_n$ that say it's first $3$ derivatives are are the standard basis vector. This means that the determinant associated to a choice of $(n_1,n_2,n_3)=(1,2,k)$ yields $k!$ times the coefficient of the $k$-th power of $\Gamma_n$. As this quantity must be $5$ also, we have that $\norm{\Gamma_n^{(k)}(0)}\lesssim 1$, for any $k$. Therefore each coefficient of $\Gamma_n$ is bounded for all $n$
, so therefore the limit, $\Gamma$, of the sequence ${\Gamma_n}_{n=1}^\infty$ does exist.

Furthermore, this limit retains the listed properties, on the set $B$ dilated about $0$ by a factor of $\frac 1 {\lim_{n\rightarrow\infty}\delta_n}$, which is the entirety of $\C$. This means that $\norm{L_\gamma(z)}\in[\frac 14,4]$ for all $z\in\C$. As $L_\gamma$ is a polynomial, it must in fact be constant. Evaluating it at $0$ yields that it must identically equal to $1$, and again, as  it is a polynomial curve, it must be given $\Gamma(z)=(z,\frac{z^2}2,\frac {z^3}6)$, up to some redistribution of the coefficients. However this clearly contradicts the fact that
$$\det\left(\Gamma_n^{(n_1)}(0),\Gamma_n^{(n_2)}(0),\Gamma_n^{(n_3)}(0)\right) =5,$$ should hold for all $n_1<n_2<n_3$, thus completing the proof, $\qed$

The second step to obtaining bounds on Dyadic Annulli involves obtaining bounds for off spring curves

\begin{lemma}
For any convex set containing $0$, $B$ there exists a decomposition of  $B=\bigcup_{j=1}^{M_{d,n}}B_j$ such that, if $\gamma$ is a degree $N$ curve satisfying

$$|L_\gamma(z)|\in[\frac12,2],\; z\in B,$$ then for $K\in\Z$, and $\bar{h}\in\C^K$, then the curve $\gamma_h(z)=\frac1K\sum_{j=1}^K\gamma(z+h_j)$ satisfies the following on $B_h=\bigcap_{j=1}^K \left(B-h_j\right)$:

$$|J_{\gamma_h}|\gtrsim \prod_{j=1}^3|L_{\gamma_h}(z_j)|^\frac13\prod_{i<j,\;j=1}^3\norm{z_i-z_j}, $$
and
$$|L{\gamma_h}(z)|\sim 1.$$
\end{lemma}

The proof of this lemma is omitted due to a lack of any noticeable change to the proof from where it appears as Lemma 2.3 in \cite{Stovall}, with the exception of the replacement of $[-1,1]$ being any convex set containing $0$, and the interval $I=[t_0-\delta,t_0+\delta]$ with balls of radius $\delta$.

Following \cite{Stovall}, we obtain our desired estimates on sets of roughly constant torsion by using a an inductive argument to get a sequence of $(p,q)$ values for which the estimate holds, and we can interpolate between these results to get the full range. The base case is the trivial estimate that $||\mathcal{E}_\gamma f||_{L^\infty}\lesssim||f||_{L^1}. $

Note that because we are on sets of roughly constant torsion, we have $||f||_{L^1(\lambda_\gamma)}\sim ||f||_{L^1}$. Our inductive step is analogous to lemma [2.4] in \cite{Stovall}.

\begin{lemma}
Suppose that, for some $p_0\in[1,7)$ we have the estimate, for some offspring curve $\gamma_h$, and set $B_h$, as described in the previous lemma, that $$||\mathcal{E}_{\gamma_h}\left(\chi_{B_h}f\right)||_{L^{6p'}}\lesssim||f||_{L^p},$$

Then we have the same estimate for $p_0$ replaced with some $p>1$ satisfying the inequality $\frac 1p>\frac {2}{15}+\frac1{p_0}$.
\end{lemma}

Note that this will obviously let us, starting from $p_0=1$, continuously get estimates for larger and larger $p$-values, with the limiting behaviour of it tending towards the desired $p=7$. Again, the proof does not change much from \cite{Stovall}, but is described here for completion.

\emph{Proof:} First, denote $v(h)=\prod_{j=1,i<j}^3(h_j-h_i)$. Furthermore, we denote the unweighted operator $$\hat{\mathcal{E}}f(\bar{z}):=\int_{B_h}e^{i\bar{z}\cdot\gamma(w)}f(w)dw.$$ We also define the measure $\mu=\mu_f$ via $$\mu(\phi):=\int_{B_h}\phi(\gamma_h(w))f(w)dw.$$ This measure is defined like this so that we have the relation $$||\hat{\mathcal{E}}_{\gamma_h}f||_{L^p}=||\hat{\mu}^3||^\frac13_{L\frac p 3}. $$
Because of this relation to the cube of the fourier transform of $\mu$, we should also consider the three-fold convolution of $\mu$ with itself, denoted $\mu^{\ast 3}(d\xi)$. From computation, we can see this self-convolution actually behaves as a function, so we denote $g(\xi)=\mu^{\ast 3}(d\xi)$, and  via computation, and the change of variable $(w_1,w_2,w_3)\rightarrow \gamma_h(w_1)+\gamma_h(w_2)+\gamma_h(w_3)$, we obtain $$\left|g(\frac13(\gamma_h(z_1)+\gamma_h(z_2)+\gamma_h(z_3)))\right|\sim\left|\frac {f(z_1)f(z_2)f(z_3)}{J_{\gamma_h}(z_1,z_2,z_3)^2}\right| .$$

Note that the square power of $2$ appears here, and not in \cite{Stovall}, due to the relation between the complex and the real Jacobian of a mapping.

For $h=(h_1,h_2,h_3)$, we adopt the notation that $h_1=0$, and define a new function $G(t,h):=g(\frac13(\gamma_h(t+h_1)+\gamma_h(t+h_2)+\gamma_h(t+h_3)))$.

Note now we have that $$\hat{g}(z)\sim\int_{\C^{2}}\int_{B_h}e^{i\bar{x}\cdot\hat{\gamma}_h(w)}g(\hat{\gamma}_h(w))|J_{\gamma_h}(t+h_1,t+\hat{h}_2,t+\hat{h}_3)|^2dwd\hat{h}_2d\hat{h}_3.$$

Here note that $\hat{\gamma}_h$ is an offspring curve of our original $\gamma_h$, generated by the vector $h=(0,\hat{h}_2,\hat{h}_3)$.

Recalling from Lemma 1.5, and as $|L_{\gamma_h}|\sim 1$ we have $|J_{\gamma_h}(z_1,z_2,z_3)|\sim|v(z_1,z_2,z_3)|$. 

This coupled with using Plancherel's Theorem with respect to $g$ gives   $$||\hat{g}||_{L^2}\lesssim||G||_{L^2_{h'}(L^2_w;|v(h)|^2)} .$$

Here, $||f||_{L^a_x(L^b_y;w(x,y))}$ simply denotes a weighted mixed norm. This quantity is what is obtained after computing the $L^b$ norm in the $y$ variable, with the weight $w(x,y)$ considered, followed by taking the $L^a$ norm in $x$ of the resulting function of $x$.

Utilising the assumption that we have the estimate at $p_0$ yields, for $q=\frac{d(d+1)}2p_0'$ $$||\hat{g}||_{L^q}\lesssim||G||_{L^1_{h'}(L^{p_0}_w;|v(h)|^2)} ,$$
and as we know we have the estimate at $q=\infty$, $p=1$
 $$||\hat{g}||_{L^\infty}\lesssim||G||_{L^1_{h'}(L^1_w;|v(h)|^2)}.$$

Using interpolation, we obtain

 $$||\hat{g}||_{L^c}\lesssim||G||_{L^a_{h'}(L^b_w;|v(h)|^2)}.$$
 for $(a,b,c)$ satisfying $$\frac5a+\frac 1b+\frac 6c=6.$$
 
 Now, we have $||G||_{L^a_{h'}(L^b_w;|v(h)|^2)}\sim\int|v(h)|^{-2(a-1)}\left(\int_{B_h}f(w+h_1)f(w+h_2)f(w+h_3)\right)$
 Now, as $\frac 1{v(h)}\in L^\frac32$, one can conclude that
$$ ||G||_{L^a_{h'}(L^b_w;|v(h)|^2)}\leq||f||_{L_t^{p,1}},$$
whenever $a\in(1,\frac 53)$, and $b\in[a,\frac{2a}{5-3a})$, such that $\frac 3p=\frac{5}{a}+\frac 1b-3$

From here, we proceed identically as in \cite{Stovall} to obtain our lemma, and therefore our main Theorem for our bounds on dyadic sets.

Now, we will deal with the case where we are not on  set where the torsion is approximately constant, and instead on a set where the torsion behaves like a monomial. In this section, we use the same decomposition of $\C$ as described in the convolution section. Recall that we have $\C=\bigcup_{j=1}^{M_N}B_j$, and on any $B_j$, we had $|L_\gamma(z)|\sim A_j|z-b_j|^{k_j}$, and $|L_1(z)|=|\gamma_1'(z)|\sim B_j|z-b_j|^{l_j}$ for some integers $l_j,\;k_j\geq 0$, real numbers $A_j,\;B_j$ and some complex number $b_j$. Also recall that we have the Jacobian bound: $$|J_{\gamma}(z_1,z_2,z_3)|\gtrsim \prod_{j=1}^3|L_{\gamma}(z_j)|^\frac13\prod_{i<j,\;j=1}^3\norm{z_i-z_j}, $$ on each set from this decomposition.

Note in this section we may assume $k_j\neq 0$, as otherwise we are on a set on which our torsion can be considered constant, which we obtained the bound for in the last section. On the sets where $k_j\neq 0$, we further subdivide them. 

Let $B_j$ be such a set, and let $B_{j,n}=\{z\in B_j|2^n<|z-b_j|<2^{n+1}\}$. Note as after an affine transformation, each $B_j$ is contained in a small sectors, the $B_{j,n}$ annuli can be ``convexified'' as in the previous section, and can be considered as triangles or trapezoids, depending if $B_j$ contained  $0$ or not.

Note this assumption relies that our affine transformation forces $b_j=0$, and that we can simultaneously ensure $A_j=B_j=1$ also. Now, we are ready to introduce the square function estimate from \cite{Stovall}.

\begin{lemma}
For each $q\in(7,\infty)$, we have the following bound, in terms of the above $B_j$ and $B_{j,n}$: $$\left|\left|\mathcal{E}_\gamma\chi_{B_j}f\right|\right|_{L^q}\lesssim \left|\left|\left(\sum_n|\mathcal{E}_\gamma\chi_{B_{j,n}}f|^2\right)^{\frac 12}\right|\right|_{L^q} $$.
\end{lemma}

\emph{Proof:} First note that if the amount of non-empty $B_{j,n}$, say $M$, is $\mathcal{O}(1)$, then the lemma immediately follows from the triangle inequality and the equivalence of $l_1$ and $l_2$ norms in $\C^M$.

We may therefore assume than $M\neq\mathcal{O}(1)$. Furthermore, if we let $a_j$ be the complex number with smallest modulus in  $B_j$, and $n_j=\lceil\log_{2}|a_j|\rceil+1$, then we can replace $B_j$ with $B_j'=B_j\bigcap\{z\in\C\;|\;|z|>2^{n_j}\}$. As this excludes at most $2$ of the $B_{j,n}$, we, as above, immediately have the estimate on these two sets.  Therefore, if we prove the case with $B_j=B_j'$, we have
$$\left|\left|\mathcal{E}_\gamma\chi_{B_j}f\right|\right|_{L^q}\leq \left|\left|\mathcal{E}_\gamma\chi_{B_j\setminus B_j'}f\right|\right|_{L^q}+\left|\left|\mathcal{E}_\gamma\chi_{B_j'}f\right|\right|_{L^q}
$$

$$\lesssim\left|\left|\left(\sum_{n=n_j-2}^{n_j-1}|\mathcal{E}_\gamma\chi_{B_{j,n}}f|^2\right)^{\frac 12}\right|\right|_{L^q} +\left|\left|\left(\sum_{n\geq n_j}|\mathcal{E}_\gamma\chi_{B_{j,n}}f|^2\right)^{\frac 12}\right|\right|_{L^q}.  $$

Now the quantity on the last line is less than the twice the maximum of these two terms, and the maximum of these two terms is less than the right hand side of our initial inequality, so proving this assumption is indeed sufficient. To see this notice that we can write $\gamma_1(z)=\int_{a_j}^z\gamma'_1(w)dw$, where we recall that this notation refers to the path integral of $\gamma_1'$ along the straight line joining $a_j$ to $z$. 

As $\gamma'_1$ is a polynomial, and our $\C$ decomposition was with respect to $\gamma_1'$, we can consider $\gamma_1'$ as a sector contained function, as in the convolution section. We therefore have $|\gamma_1(z)|\sim|z|^{l_j+1}-|a_j|^{l_j+1}\sim|z|^{l_j+1}$, due to the separation ensured by restricting ourselves to $B_j'$.

Now, we make the claim that has its support contained within $\text{Im}\left(\gamma\right)\times \C^2$. That is to say it is contained within $\left\{\xi\in\C^3||\xi_1|\in[2^{nl_j},2^{nl_j+1}]\right\}$. To see this, consider $\langle\widehat{\mathcal{E}_\gamma\chi_{B_{n,j}}f},\phi\rangle$ where $\phi$ is a Schwarz functions with support contained entirely outside of the image of $\gamma$. We have:

\begin{equation*}
\begin{split}
\langle\widehat{\mathcal{E}_\gamma\chi_{B_{n,j}}f},\phi\rangle&=\langle\mathcal{E}_\gamma\chi_{B_{n,j}}f,\widehat{\phi}\rangle \\
 & =\int_\C \mathcal{E}_\gamma\chi_{B_{n,j}}f(z)\widehat{\phi}(z)\text{dz}\\
 &=\int_\C \int_{B_{n,j}} e^{i(z\cdot\gamma(w)}f(w)\lambda_\gamma(w)\text{dw}\widehat{\phi}(z)\text{dz} \\
 &=\int_{B_{n,j}}f(w)\lambda_\gamma(w) \int_\C e^{i(z\cdot\gamma(w))}\widehat{\phi}(z)\text{dz}\text{dw}\\
 &=\int_{B_{n,j}}f(w)\lambda_\gamma(w) \phi(\gamma(w))\text{dw},
\end{split}
\end{equation*}

by the Fourier inversion theorem. As we chose $\phi$ to not have support in the image of $\gamma$, this evaluate to zero. As this holds for all $\phi$ satisfying this condition, this means the support of $\widehat{\mathcal{E}_\gamma\chi_{B_{n,j}}f}$ is contained in the image of $\gamma$. Therefore we can use the following theorem which appears in \cite{Stein}. Here, the operator $S_\rho$ denotes the operator such that $\widehat{S_\rho f}=\chi_\rho\hat f$, and $\Delta$ is the collection of rectangles from the dyadic decomposition of $\C^3$, or equivalently $\R^{6}$, which are just rectangles with sides which are products of 6 intervals of the form $[2^{k_i},2^{k_i+1}]$ for some various integer $k$ values . Theorem 5. in Chapter $IV$ of \cite{Stein} is restated here.

\begin{thm}
If $f\in L^p(\R^n)$, then $$\left(\sum_{\rho\in\Delta} |S_\rho f|^2\right)^\frac 12\in L^p\left(\R^n\right) \text{ and }||f||_{L^p}\sim\left|\left|\left(\sum_{\rho\in\Delta} |S_\rho f|^2\right)^\frac 12\right|\right|_{L^p}.$$
\end{thm}

Notice now our lemma immediately follows from letting $\mathcal{E}_\gamma\chi_{B_{j,n}}f$ serve the role of $f$ in this theorem, and further notice that $S_\rho \mathcal{E}_\gamma\chi_{B_{j}}f=0$ if the rectangle $\rho$ does not have it's edge parallel to the first axis project onto the interval $[2^{nl_j},2^{nl_j+1}]$ for some $n$ corresponding to a non-empty $B_j,n$.

Furthermore, if $\rho$ does satisfy the above stated condition, we have $S_\rho \mathcal{E}_\gamma\chi_{B_{j}}f=\mathcal{E}_\gamma\chi_{B_{j,n}}f$, for the $n$ value as stated in the condition, which proves the lemma.$\qed$

\emph{Proof of Theorem 1.3:} To prove Theorem 1.3, from here, we first restrict ourselves to considering $q<12$. If we establish the theorem for $q<12$, we can interpolate between any estimate corresponding to some arbitrary $q<12$ and the trivial estimate at $q=\infty$ to get the full range. 

Note furthermore that by the triangle inequality, we need only establish the estimate for $f$ being a function which is only supported in a single $B_j$. As this $j$ is now fixed, we change our notation to label $B_j=B$, and $B_{j,n}=B_n$

We now have, by our previous lemma:

$$||\mathcal{E}_\gamma f||^q_{L^q}\lesssim\int\left(\sum_n|\mathcal{E}_\gamma \chi_{B_n}f(x)|^2\right)^\frac q2dx=\int\prod_{j=1}^6\left(\sum_{n_j}|\mathcal{E}_\gamma \chi_{B_{n_j}}f(x)|^2\right)^\frac q{12}dx$$

$$\lesssim\int\prod_{j=1}^6\sum_{n_j}|\mathcal{E}_\gamma \chi_{B_{n_j}}f(x)|^\frac q6dx\sim\sum_{n_1<n_2\cdots<n_6}\int\prod_{j=1}^6|\mathcal{E}_\gamma \chi_{B_{n_j}}f(x)|^\frac q6dx.$$

Next, we introduce a lemma whose proof will be postponed until the main theorem has been established:

\begin{lemma}
There exists a real, positive number $\varepsilon$ depending only on $N$ and $p$ such that under the assumptions made of $\gamma$, $B$, $q$ and $p$. Then for $n_1<n_2\cdots<n_6$, and $f_j\in L^p$ are supported on $n_j$, we have that: $$||\prod_{j=1}^6\mathcal{E}_\gamma f_j||_{L^q}^q\lesssim 2^{-\varepsilon(n_6-n_1)}\prod_{j=1}^6||f_j||_{L^p(\lambda_\gamma)}. $$
\end{lemma}

Combining this lemma with the preceding inequality yields

$$||\mathcal{E}_\gamma f||^q_{L^q}\lesssim\sum_{n_1<n_2\cdots<n_6}2^{-\varepsilon(n_6-n_1)}\prod_{j=1}^6||\chi_{B_{n_j}}f||^{\frac q6}_{L^p(\lambda_\gamma)}.$$

Let us denote $B_{a,b}=\{z\in B\hspace{.1cm}|\hspace{.1cm}|z|\in [2^a,2^b]\}$. Of course as a consequence, $B_{n_j}\subset B_{n_1,n_6}$. Using this, letting $m=n_6-n_1$, and $n=n_1$, and then noticing that as there are a bounded number of choices for $n_2,\cdots,n_5$ for fixed $m$, we can rewrite our estimate as, for some constant $C$, we have:
$$||\mathcal{E}_\gamma f||^q_{L^q}\lesssim\sum_m\sum_n2^{-\varepsilon m}m^{C}||\chi_{ B_{n,n+m}}f||^q_{L^p(\lambda_\gamma)}.$$

Now, for any fixed $m$, we have:
$$\sum_n||\chi _{B_{n,n+m}}f||^q_{L^p(\lambda_\gamma)}\lesssim\left(\sup_n||\chi_{B_{n,n+m}}f||_{L_p(\lambda_\gamma)}^{q-p}\right)\sum_n||\chi_{ B_{n,n+m}}f||^q_{L^p(\lambda_\gamma)}.$$
Now, as any point can only be in at most $m$ sets of the form $B_{n,n+m}$, we can observe:
$$\sum_n||\chi _{B_{n,n+m}}f||^q_{L^p(\lambda_\gamma)}\lesssim m||f||_{L^p(\lambda_\gamma)}^q. $$

Combining this with our main chain of inequalities leads to: $$||\mathcal{E}_\gamma f||^q_{L^q}\lesssim \sum_m 2^{-\varepsilon m}m^{C+1}||f||^q_{L^p(\lambda_\gamma)}\sim||f||^q_{L^p(\lambda_\gamma)}.$$ 

So we are done $\hfill\qed$.

Now, to return to the proof of the postponed lemma.\\
\emph{Proof of Lemma 6.1:} By Hölder’s inequality, we have that:$$||\prod_{j=1}^6\mathcal{E}_\gamma f_j||_{L^{\frac q6}}\leq \left(\prod_{j=4}^6||\mathcal{E}_\gamma f_{i_j} ||_{L^q}\right)\left(||\prod_{j=1}^3\mathcal{E}_\gamma f_{i_j} ||_{L^\frac q3}\right), $$

where the $i_j$  values are some enumeration of $\{1,2\cdots 6\}$. Notice now that the lemma will be proven if we can establish that: $$||\prod_{j=1}^3\mathcal{E}_\gamma f_{j} ||_{L^\frac q3}\lesssim 2^{-\varepsilon(n_6-n_1)}\prod_{j=1}^3||f_j||_{L^p(\lambda_\gamma)}
. $$

Now, notice that we have, by Hölder’s inequality, and by the fact that the support of each $f_j$ is in B, and Theorem 4.3, that:
$$||\prod_{j=1}^3\mathcal{E}_\gamma f_{j} ||_{L^\frac q3}\leq\prod_{j=1}^3||\mathcal{E}_\gamma f_{j} ||_{L^q}\lesssim \prod_{j=1}^3||f_j||_{L^p(\lambda_\gamma)}.$$

Therefore the lemma  follows immediately if $n_3-n_1$ is bounded by some arbitrary constant. That is to say, if $n_6-n_1\leq 6$ we are done immediately. Therefore, we can assume $n_3-n_1\geq 6$. 

Furthermore, we can prove this for $q=12$, and $p=2$, as if this is established we can interpolate with between this $(p,q)$ value and some $(p,q)$ value arbitrarily close to the endpoint. So the problem is now to prove that: $$||\prod_{j=1}^3\mathcal{E}_\gamma f_{j} ||_{L^4}\lesssim 2^{-\varepsilon(n_6-n_1)}\prod_{j=1}^3||f_j||_{L^2(\lambda_\gamma)}.
 $$
 
 Now, we define the measure $d\mu_j$ such that: $$d\mu_j(\phi):=\int_{B_{n_j}}\phi(\gamma(z))f_j(z)\lambda_\gamma(z)dz.$$
 
 Then, by the Haussdorff-Young theorem, we have: $$||\prod_{j=1}^3\mathcal{E}_\gamma f_{j} ||_{L^4}\leq||d\mu_1\ast d\mu_2\ast d\mu_3 ||_{L^{\frac 43}}.$$
 
 Now, we can compute this 3-fold convolution explicitly as:
 $$d\mu_1\ast d\mu_2\ast d\mu_3(\phi)=\int_{B_{n_1}}\int_{B_{n_2}}\int_{B_{n_3}}\phi\left(\sum_{j=1}^3\gamma(z_j)\right)\left(\prod_{j=1}^3f_j(z_j)\lambda_\gamma(z_j)\right)dt_3dt_2dt_1.$$
 
We now rewrite this integral in terms of the permutation group $S_3$, and the sets \newline$P_\sigma=\{z\in B_{n_1}\times B_{n_2}\times B_{n_3}\hspace{.1cm}|\hspace{.1cm}|z_{\sigma(1)}|<|z_{\sigma(2)}|<|z_{\sigma(3)}|\}$, on which the function \newline$\Phi_\gamma(z_1,z_2,z_3)=\sum_{j=1}^3\gamma(z_j)$ is injective. So we have:

 $$d\mu_1\ast d\mu_2\ast d\mu_3(\phi)=\sum_{\sigma\in S_3}\int_{P_\sigma}\phi\left(\sum_{j=1}^3\gamma(z_j)\right)\left(\prod_{j=1}^3f_j(z_j)\lambda_\gamma(z_j)\right).$$

 Therefore, by making the substitution $\nu=\phi_\gamma(z_1,z_2,z_3)$ we can notice that this 3-fold convolution of measure is actually a function of $\C^3$ given by:
 
$$d\mu_1\ast d\mu_2\ast d\mu_3(\nu)=\sum_{\sigma\in S_3}F_\sigma(\nu), $$

where $$F_\sigma(\nu):=\chi_{P_\sigma}(z)\left(\prod_{j=1}^3f_j(z_j)\lambda_\gamma(z_j)\right)|J_\gamma(z)|^{-2}. $$

Note here that we proceed with understanding that $z$ is just the preimage of $\nu$ under the mapping $\Phi_\gamma$ restricted to $P_\sigma$. We therefore have:
$$||F_\sigma||_{L^\frac 43}=\left(\int_{P_\sigma}\left(\prod_{j=1}^3|f_j(z_j(\nu))\lambda_\gamma(z_j(\nu))||J_\gamma(z(\nu))|^{-2}\right)^\frac 43d\nu \right)^\frac 34$$
$$ =\left(\int_{P_\sigma}\left(|\prod_{j=1}^3f_j(z_j)\lambda_\gamma(z_j)||J_\gamma(z)|^{-2}\right)^\frac 43 |J_\gamma(z)|^{2}dz \right)^\frac 34$$

$$=\left(\int_{P_\sigma}\left(\prod_{j=1}^3|f_j(z_j)\lambda_\gamma(z_j)||J_\gamma(z)|^{-\frac 12}\right)^\frac 43 dz \right)^\frac 34$$
$$=||\chi_{P_\sigma}\prod_{j=1}^3f_j(z_j)\lambda_\gamma(z_j)J_\gamma(z)^{-\frac 12}||_{L^{\frac 43}}$$
$$\lesssim||\chi_{P_\sigma}\prod_{j=1}^3f_j(z_j)\lambda_\gamma(z_j)^\frac 12\prod_{1\leq i<j\leq3}|z_i-z_j|^{-\frac 12}||_{L^{\frac 43}}, $$
from applying our Jacobian bound. Now by the pigeonhole principle, some consecutive $n_j$ and $n_{j+1}$ are closer together than the average difference between the $n_k$ values. That is to say, for $k=1$ or $k=2$ we have $n_{k+1}-n_k\geq\frac{n_3-n_1}3\geq 2$. We can do both cases simultaneously, understanding an empty product represents $1$.

We have $$\prod^3_{j=1,\hspace{.1cm}i< j}|z_j-z_i|=\prod_{i\leq k\hspace{.1cm}j\geq k+1}|z_j-z_i|\prod^k_{j=1,\hspace{.1cm}i< j}|z_j-z_i|\prod^3_{i=k,\hspace{.1cm}i< j}|z_j-z_i|$$
$$\sim\prod_{i\leq k\hspace{.1cm}j\geq k+1}2^{n_j}\prod^k_{j=1,\hspace{.1cm}i< j}|z_j-z_i|\prod^3_{i=k,\hspace{.1cm}i< j}|z_j-z_i| . $$ 

This yields:
$$||F_\sigma||^{\frac 43}_{L^\frac 43}= 2^{-\frac k6(\sum_{j=k}^2 n_{k+1})}T_k(f_1,\cdots f_k)T_{3-k}(f_{k+1},\cdots f_3)$$

Where we have defined $$T_l(f_1,\cdots f_l):=\int_{\C^l}\prod_{j=1}^lf_j(z_j)^\frac 43\lambda_\gamma(z_j)^\frac 23\prod_{j=1,\hspace{.1cm}i<j}^l |z_j-z_i|^{-2}dz$$.

Note, in \cite{Christ2}, Christ establishes the fact that $$\int_{\R^l}\prod_{i=1}^lf_i(t_i)\prod_{1\leq i<j\leq l}g_{i,j}(t_i-t_j)\lesssim \prod_{i=1}^l||f_i||_{L^p}\prod_{1\leq i<j\leq l}||g_{i,j}||_{L^{q,\infty}}$$

for $p\in[1,l)$, and $q$ satisfying $p^{-1}+\frac{l-1}{2}q^{-1}=1$. The translation of this argument to the Complex Case can be done without any significant changes to the argument presented there. Furthermore,  if $l=1$, the desired inequality is trivial. Applying this to $T_k$ with $p=\frac{6}{7-k}$, which yields $q=3$ for each k, we get:

$$T_k(z_1,\cdots,z_k)\lesssim \prod_{i=1}^3||f_i\lambda_\gamma^{\frac 12}||^\frac{4}{3}_{L^{\frac{6}{7-k}}}$$

Recalling that the support of $f_j$ is contained within $B_{n_j}$, and we have $|B_{n_j}|\sim2^{n_j}$. Utilising this with Hölder's inequality gives us:$$T_k\lesssim\prod_{j=1}^k 2^\frac{n_j(3-k)}{6}||f_j||^\frac{4}{3}_{L^2(\lambda_\gamma)} $$
and $$T_{3-k}\lesssim\prod_{j=k+1}^3 2^\frac{n_jk}{6}||f_j||^\frac{4}{3}_{L^2(\lambda_\gamma)} .$$

Combining this with our bound for the $L^\frac 43$ norm of $F_\sigma$ yields that: $$||F_\sigma||_{L^{\frac 43}}^\frac 43\lesssim 2^{\frac{(n_1+\cdots+n_k)(3-k)}{6}-\frac{(n_{k+1}+\cdots +n_3)k}{6}}\prod_{i=1}^3||f_i||_{L^2(\lambda_\gamma)}$$
$$\lesssim 2^{-\frac{k(3-k)}{6}(n_{k}-n_{k+1})}\prod_{i=1}^3||f_i||_{L^2(\lambda_\gamma)}.$$

Therefore, as $n_{k+1}-n_k\geq \frac{n_3-n_1}{3}\geq2$, this completes the proof of the lemma. \qed

\end{document}